\documentclass[11pt,reqno]{article}
\usepackage{latexsym, amsmath, amssymb, a4, epsfig}
\usepackage{graphicx}

\newtheorem{theorem}{Theorem}[section]
\newtheorem{cor}[theorem]{Corollary}
\newtheorem{lemma}[theorem]{Lemma}

\newtheorem{example}{Example}

\newcommand{\nm}{\noalign{\smallskip}}
\newcommand{\qed}{ $\Box$}
\newcommand{\ds}{\displaystyle}
\newcommand{\pf}{\noindent {\sl Proof}. \ }
\newcommand{\p}{\partial}
\newcommand{\pd}[2]{\frac {\p #1}{\p #2}}
\newcommand{\eqnref}[1]{(\ref {#1})}

\newcommand{\Cbb}{\mathbb{C}}
\newcommand{\Dbb}{\mathbb{D}}
\newcommand{\Ebb}{\mathbb{E}}
\newcommand{\Fbb}{\mathbb{F}}

\newcommand{\Kbb}{\mathbb{K}}
\newcommand{\Rbb}{\mathbb{R}}
\newcommand{\Sbb}{\mathbb{S}}

\newcommand{\la}{\langle}
\newcommand{\ra}{\rangle}

\newcommand{\Ecal}{\mathcal{E}}
\newcommand{\Hcal}{\mathcal{H}}

\newcommand{\Kcal}{\mathcal{K}}
\newcommand{\Ncal}{\mathcal{N}}
\newcommand{\Dcal}{\mathcal{D}}

\newcommand{\Rcal}{\mathcal{R}}
\newcommand{\Scal}{\mathcal{S}}

\def\Bk{{\bf k}}

\def\Bs{{\bf s}}

\def\Bx{{\bf x}}

\def\BH{{\bf H}}
\def\BI{{\bf I}}

\newcommand{\Ga}{\alpha}
\newcommand{\Gb}{\beta}
\newcommand{\Gd}{\delta}
\newcommand{\Ge}{\epsilon}

\newcommand{\Gvf}{\varphi}
\newcommand{\Gg}{\gamma}

\newcommand{\Gl}{\lambda}

\newcommand{\Gv}{\nu}

\newcommand{\Gt}{\theta}

\newcommand{\Gs}{\sigma}

\newcommand{\Go}{\omega}

\newcommand{\GD}{\Delta}

\newcommand{\GG}{\Gamma}
\newcommand{\GL}{\Lambda}

\newcommand{\GT}{\Theta}
\newcommand{\GS}{\Sigma}

\newcommand{\GO}{\Omega}

\newcommand{\beq}{\begin{equation}}
\newcommand{\eeq}{\end{equation}}

\numberwithin{equation}{section}
\numberwithin{figure}{section}

\begin{document}
\title{Spectral properties of the Neumann-Poincar\'e operator and uniformity of estimates for the conductivity equation with complex coefficients\thanks{\footnotesize This work is supported by the Korean Ministry of Education, Sciences and Technology through NRF grants Nos. 2010-0017532 and 2013R1A1A1A05009699.}}

\author{Hyeonbae Kang\thanks{Department of Mathematics, Inha University, Incheon
402-751, S. Korea (hbkang, kskim, hdlee@inha.ac.kr).} \and Kyoungsun Kim\footnotemark[2] \and Hyundae Lee\footnotemark[2] \and Jaemin Shin\thanks{Department of Mathematical Sciences, Hanbat National University, Daejeon 305-719, S. Korea (jaemin.shin@hanbat.ac.kr).} \and Sanghyeon Yu\thanks{Seminar for Applied Mathematics, ETH Z\"urich,
R\"amistrasse 101, CH-8092 Z\"urich, Switzerland (sanghyeon.yu@sam.math.ethz.ch).} }

\maketitle

\begin{abstract}
We consider well-posedness of the boundary value problem in presence of an inclusion with complex conductivity $k$. We first consider the transmission problem in $\Rbb^d$ and characterize solvability of the problem in terms of the spectrum of the Neumann-Poincar\'e operator. We then deal with the boundary value problem and show that the solution is bounded in its $H^1$-norm uniformly in $k$ as long as $k$ is at some distance from a closed interval in the negative real axis. We then show with an estimate that the solution depends on $k$ in its $H^1$-norm Lipschitz continuously. We finally show that the boundary perturbation formula in presence of a diametrically small inclusion is valid uniformly in $k$ away from the closed interval mentioned before. The results for the single inclusion case are extended to the case when there are multiple inclusions with different complex conductivities: We first obtain a complete characterization of solvability when inclusions consist of two disjoint disks and then prove solvability and uniform estimates when imaginary parts of conductivities have the same signs. The results are obtained using the spectral property of the associated Neumann-Poincar\'e operator and the spectral resolution.
\end{abstract}

\noindent{\footnotesize {\bf AMS subject classifications}. 35J47 (primary), 35P15 (secondary)}

\noindent{\footnotesize {\bf Key words}. Neumann-Poincar\'e operator, Lipschitz domain, spectrum, spectral resolution, asymptotic formula, uniformity, complex conductivity, high contrast}

\section{Introduction}\label{sec:intro}
Let $\GO$ be a bounded simply connected domain in $\Rbb^d$ ($d\geq 2$) and let $D$ be a simply connected domain compactly contained in $\GO$. We assume that boundaries of $D$ and $\GO$ are Lipschitz continuous and the conductivity of $\GO \setminus D$ is $1$ while that of $D$ is $k$ so that the conductivity distribution is given by
\beq
\Gg_k = \chi(\GO \setminus D)+ k \chi(D),
\eeq
where $\chi(D)$ is the characteristic function of $D$. For a given Neumann data $g$ consider the following elliptic problem:
\beq\label{NBP}
 \left \{
 \begin{array}{l}
 \ds \nabla  \cdot ( \Gg_k  \nabla u) =0 \quad \mbox{in } \GO,  \\
 \nm \ds \p_\nu u |_{\p\GO} =g, \\
 \nm \ds \int_{\p\GO} u \, d\Gs=0.
 \end{array}
 \right.
 \eeq
Here and throughout this paper $\p_\nu u$ denotes the outward normal derivative of $u$ and the conductivity $k=k'+ik''$ is a complex number ($k'$ and $k''$ denote the real and imaginary parts of $k$, respectively). We emphasize that the problem \eqnref{NBP} admits a unique solution if $k$ is not on the negative real axis.

The solution $u_k$ to \eqnref{NBP} varies depending on the conductivity $k$ regarded as a parameter. For instance, we have a standard regularity estimate
\beq\label{standard}
\| u_k \|_{H^1(\GO)} \le C \| g \|_{H^{-1/2}(\p\GO)}
\eeq
for some constant $C$ which may depend on $k$. Here and throughout this paper $H^s(\GO)$ (and $H^s(\p \GO)$) denotes the standard $L^2$-Sobolev space and $H^{-1/2}(\p\GO)$ is the dual space of $H^{1/2}(\p\GO)$. However, if $k$ is real ($0 \le k \le \infty$), it is proved in \cite{NgVo09} that the constant $C$ can be chosen independently of $k$. Moreover, if the inclusion $D$ is diametrically small, it is proved in the same paper that the asymptotic boundary perturbation formula (see section \ref{sec:asymp}) is valid uniformly in $k$.
We emphasize that these results were obtained using variational methods.

\begin{figure}[h!]
\begin{center}
\epsfig{figure=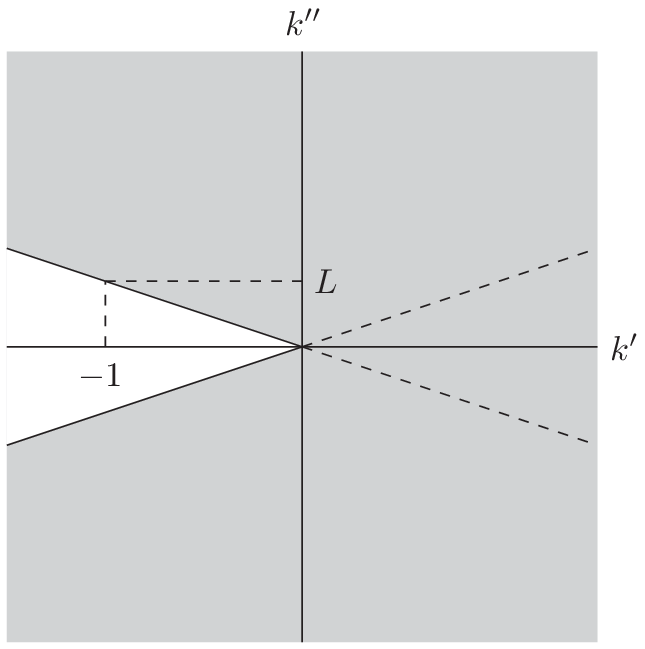,height=5cm, width=5cm}\hskip 1cm
\epsfig{figure=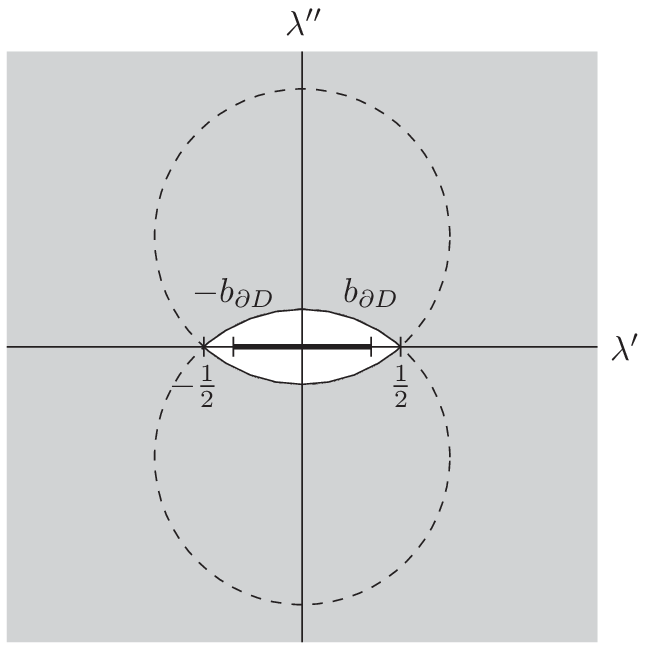,height=5cm, width=5cm}
\end{center}
\caption{Left: the region $|k''| \ge - L k'$ (grey). Right: the region after the transformation $\Gl=\frac{k+1}{2(k-1)}$. $b_{\p D}$ is the spectral bound of the NP operator.}\label{fig}
\end{figure}

In this paper we develop a new method to investigate solvability of \eqnref{NBP} and dependency of solution $u_k$ on the conductivity $k$ of the inclusion, and to extend above mentioned results to the case when $k$ is a complex number. We first prove existence and uniqueness of the solution to \eqnref{NBP} when $k'' \neq 0$. We then show that if $k=k' + i k''$ satisfies
\beq\label{kcondition}
|k''| \ge - L k'
\eeq
for any given constant $L>0$ (see the left figure in Figure \ref{fig}), we prove that \eqnref{standard} holds for some $C$ independent of $k$. We then show \eqnref{standard} holds near $k=0$ and $k=\infty$ by a perturbation argument. The uniform estimates near $0$ and $\infty$, and that in the sector imply that there is a closed interval such that uniform estimate holds away from the interval. We also show with a precise quantitative estimate that $u_k$ depends on $k$ in $H^1(\GO)$-norm Lipschitz continuously. We finally prove uniform validity of the asymptotic boundary perturbation formula when $D$ is diametrically small, regardless of $k$ as long as $k$ satisfies \eqnref{kcondition}.

Solvability and uniformity results for the single inclusion case are extended to the case when there are multiple inclusions. Namely, there are multiple inclusions $D_j$, $j=1, \ldots, M$, with complex conductivities $k_j$. If all the conductivities are the same, namely, $k_1=\ldots=k_M$, then results for the single inclusion case are valid for the multiple inclusion case without change. However, if we allow them to be different (and complex), finding meaningful conditions which guarantee solvability of the problem seems quite difficult (see subsection \ref{sec:twodisk}). We attempt to present one sufficient condition in this paper. The condition is basically that $k_j''$ have the same sign for all $j$.

There is growing interest in the complex conductivity, especially in relation to the electrical impedance tomography (EIT). For instance, the imaginary part of the complex conductivity changes depending on the frequency of the prescribed current, and by exploiting this property inclusion can be reconstructed with a high resolution (see, for example, \cite{kimseo}). Another example is the size estimation problem in the EIT, which is to derive bounds on the volume fraction of the inclusion via boundary measurements. This problem is recently considered when the conductivity of the inclusion is complex \cite{BFV11} \cite{KKLLM14} \cite{TM13}. The bounds for the complex conductivity case turned out to be quite tight as the numerical examples presented in the last two papers show.

The method of this paper is based on the spectral property of the Neumann-Poincar\'e (NP) operator related to the problem \eqnref{NBP}.
The NP operator is a boundary integral operator which appears naturally when solving the Neumann (and Dirichlet) boundary value problems using single or double layer potentials. It is not self-adjoint with respect to the usual $L^2$-inner product. However, it can be symmetrized using a certain twisted inner product and Plemelj's symmetrization principle \cite{KhPuSh07}. We show that if $\p D$ is Lipschitz, then the spectrum of the NP operator defined on $\p D$ on a Sobolev space lies in $[-b_{\p D}, b_{\p D}]$ for some $b_{\p D} < 1/2$. Here, the spectrum or resolvent $\Gl$ of the NP operator is related to the conductivity $k$ by the bilinear transformation
\beq\label{bilinear}
\Gl=\Gl(k)=\frac{k+1}{2(k-1)} .
\eeq
Since bilinear transformations map circles on the Riemann sphere onto circles, this transformation maps the region defined by \eqnref{kcondition} onto the region outside the oval shaped curve (see Figure \ref{fig}). We emphasize that the boundary of the region in the $\Gl$-space intersects with the real axis at $1/2$ and $-1/2$. Since $b_{\p D} < 1/2$, there is some distance between the transformed region and the spectrum of the NP operator. Using this property and the spectral resolution of a self-adjoint operator we are able  to obtain the results described above.

We also consider the following transmission problem in $\Rbb^d$:
\beq\label{trans}
 \left \{
 \begin{array}{l}
 \ds \nabla  \cdot ( \chi(\Rbb^d \setminus D) + k \chi(D))   \nabla u =0 \quad \mbox{in } \Rbb^d,  \\
 \nm u(x)-h(x) = O(|x|^{1-d}) \quad\mbox{as } |x| \to \infty,
 \end{array}
 \right.
 \eeq
where $h$ is a given harmonic function in $\Rbb^d$. This problem is simpler than \eqnref{NBP} because the function $h$ does not depend on $k$ (see the next section), and solvability of this problem is completely characterized by the spectrum of the NP operator. In fact, we prove that the problem is well-posed if $\Gl(k)$ does not belong to the spectrum of the NP operator (see section \ref{sec:trans}). We emphasize that the condition \eqnref{kcondition} is required for the boundary value problem to show that the Neumann-to-Dirichlet map for the problem \eqnref{NBP} is bounded uniformly in $k$ (see Lemma \ref{Di-N}).

The program of estimating the spectrum of the NP operator goes back to Poincar\'e \cite{Po96} as the name suggests, and was recently revisited in \cite{KhPuSh07} with a new perspective. The spectrum of the NP operator also has been studied using complex function theory such as the Beurling-Ahlfors transform and quasi-conformal mapping,  for which we refer to \cite{PePu12} (and references therein) where a bound on the essential spectrum of the NP operator for two dimensional domains with corners is obtained. It is proved in \cite{Ve84} that the set $\mathbb{R}\setminus (-1/2,1/2)$ is contained in the resolvent set of the NP operator on $L^2_0(\p\GO)$ for a domain $\GO$ with the Lipschtz boundary. Its spectrum on $H^{-1/2}$ lies in $(-1/2, 1/2]$ (see \cite{ChLe08}). Recently the spectral theory of the NP operator has been applied to analysis of cloaking by anomalous localized resonance on the plasmonic structure \cite{ACKLM}, analysis of high concentration of the gradient \cite{ACKLY} \cite{BT} \cite{BT2} \cite{LiYu14}, and a shape optimization problem \cite{ACLZ}.

This paper is organized as follows. In section \ref{sec:spec} we review symmetrization and spectral properties of the Neumann-Poincar\'e operator and prove equivalence of norms. In section \ref{sec:trans} we deal with the transmission problem in the whole space $\Rbb^d$. In section \ref{sec:bvp} we deal with the boundary value problem and show that the solution exists and is unique, it is bounded in its $H^1$-norm uniformly in $k$ away from a closed interval in the negative real axis, and a quantitative estimate is derived which shows that the solution depends Lipschitz continuously on $k$. We then show in section \ref{sec:asymp} that the boundary perturbation formula when the inclusion is diametrically small is valid uniformly in $k$. The last section is to deal with the case when there are multiple inclusions.

\section{Spectral properties of the Neumann-Poincar\'e operator}\label{sec:spec}

Let $\GG(x)$ be the fundamental solution to the
Laplacian, {\it i.e.},
\beq\label{gammacond}
\GG (x) =
\begin{cases}
\ds \frac{1}{2\pi} \ln |x|\;, \quad & d=2 \;, \\ \nm \ds
\frac{1}{(2-d)\Go_d} |x|^{2-d}\;, \quad & d \geq 3 \;,
\end{cases}
\eeq
where $\Go_d$ denotes the area of the unit sphere in $\Rbb^d$. Suppose that the inclusion $D$ has a single simply connected component. The single layer potential $\Scal_{\p D} [\Gvf]$ of a density function $\Gvf \in H^{-1/2}(\p D)$ is defined by
\beq
\Scal_{\p D} [\Gvf] (x) := \int_{\p D} \GG (x-y) \Gvf (y) \, d\Gs(y) \;, \quad x \in \Rbb^d .
\eeq
It is well known (see for example \cite{AmKa07Book2} \cite{Fo95}) that $\Scal_{\p D} [\Gvf]$ satisfies the jump relation
\beq\label{singlejump}
\frac{\p}{\p\nu} \Scal_{\p D} [\Gvf] \Big|_\pm (x) = \biggl( \pm \frac{1}{2} I + \Kcal_{\p D}^* \biggr) [\Gvf] (x),
\quad x \in \p D\;,
\eeq
where $\frac{\p}{\p\nu}$ denotes the outward normal derivative, the subscripts $\pm$ indicate the limit from outside and inside $D$, respectively, and the operator $\Kcal_{\p D}^*$ is defined by
\beq\label{introkd}
\Kcal_{\p D}^* [\Gvf] (x) =
\frac{1}{\Go_d} \int_{\p D} \frac{(x-y)\cdot\nu_x }{|x-y|^d} \Gvf(y)\,d\Gs(y)\;, \quad x \in \p  D.
\eeq
Here $\nu_x$ denotes the outward unit normal vector to $\p D$ at $x$.

The operator $\Kcal_{\p D}^*$ is called the {\it Neumann-Poincar\'e} (NP) operator associated with the domain $D$. If $\p D$ is Lipschitz continuous, then $\Kcal_{\p D}^*$ is a singular integral operator and known to be bounded on $L^2(\p D)$ (and on $H^{-1/2}(\p D)$) \cite{CMM82}. Let $H^{-1/2}_0(\p D)$ be the collection of $\Gvf \in H^{-1/2}(\p D)$ with the mean zero, {\it i.e.},
\beq
\la \Gvf, 1 \ra =0
\eeq
where $\la \ , \ \ra$ is the $H^{-1/2}$-$H^{1/2}$ pairing. For $\Gvf, \psi \in H^{-1/2}_0(\p D)$,  define
\beq\label{innerp}
\la \Gvf,\psi \ra_\Hcal:= -\la \Gvf,  \Scal_{\p D}[\psi] \ra .
\eeq
Since $\Scal_{\p D}$ maps $H^{-1/2}(\p D)$ into $H^{1/2}(\p D)$, the right hand side of \eqnref{innerp} is well-defined. It is known (see, for example, \cite{ACKLM} \cite{BT2} \cite{Kang} \cite{KhPuSh07}) that $\la \ , \ \ra_\Hcal$ is an inner product on $H^{-1/2}_0(\p D)$ and $\Kcal^*_{\p D}$ is self-adjoint with respect to this inner product, which is a consequence of Plemelj's symmetrization principle (also known as Calder\'on's identity)
\beq\label{Plemelj}
\Scal_{\p D} \Kcal^*_{\p D} = \Kcal_{\p D} \Scal_{\p D}.
\eeq
Here $\Kcal_{\p D}$ is the adjoint of $\Kcal^*_{\p D}$ with respect to the usual $L^2$-inner product.
It is worth mentioning that $\Kcal^*_{\p D}$ is not self-adjoint in the usual inner product unless the domain $D$ is a disk or a ball \cite{Li01}. Let $\Hcal=\Hcal(\p D)$ be the space $H^{-1/2}_0(\p D)$ equipped with the inner product $\la \ , \ \ra_\Hcal$. We denote the norm associated with
$\la \ , \ \ra_\Hcal$ by $\| \, \cdot\,\|_\Hcal$.

Let us write $\| \Gvf \|_{H^s(\p D)}$ as $\| \Gvf \|_{s}$ from now on for ease of notation.

\begin{lemma}\label{lemA}
Let $D$ be a bounded Lipschitz domain in $\Rbb^d$ ($d \ge 2$). The operator $A: H^{-1/2}(\p D) \times \Cbb \to H^{1/2}(\p D) \times \Cbb$ by \beq\label{defA}
A(\psi, a) := \left( \Scal_{\p D}[\psi] + a, \la \psi, 1 \ra \right)
\eeq
is invertible.
\end{lemma}
\pf In two dimensions, it is shown in \cite[Theorem 2.26]{AmKa07Book2} that $A: L^2(\p D) \times \Cbb \to H^{1}(\p D) \times \Cbb$ is invertible. If $(\psi,a) \in L^2(\p D) \times \Cbb$ and $(\Gvf, b) \in H^{-1}(\p D) \times \Cbb$, then
\begin{align*}
\la (\Gvf, b), A(\psi, a) \ra &= \la \Gvf, \Scal_{\p D}[\psi] + a \ra + b \int_{\p D} \overline{\psi} d\Gs \\
&= \la \Scal_{\p D}[\Gvf], \psi \ra + \la \Gvf, a \ra + b \int_{\p D} \overline{\psi} d\Gs \\
&= \la \Scal_{\p D}[\Gvf]+b, \psi \ra + \overline{a} \la \Gvf, 1 \ra  = \la A(\Gvf, b), (\psi, a) \ra.
\end{align*}
So, by duality, $A: H^{-1}(\p D) \times \Cbb \to L^2(\p D) \times \Cbb$ is invertible. By interpolation we infer that $A: H^{-1/2}(\p D) \times \Cbb \to H^{1/2}(\p D) \times \Cbb$ in invertible.

In three or higher dimensions, it is proved in \cite{Ve84} that $\Scal_{\p D}: L^2(\p D) \to H^{1}(\p D)$ is invertible. So, $A$ is a Fredholm operator of index zero. Thus one can show by the exactly the same proof as in two dimensions that $A: L^2(\p D) \times \Cbb \to H^{1}(\p D) \times \Cbb$ is invertible in three dimensions. So we obtain the desired result by the same argument (duality and interpolation). \qed

Lemma \ref{lemA} says that for any $f \in H^{1/2}(\p D)$ and $b \in \Cbb$ there is a unique pair $(\psi,a) \in H^{-1/2}(\p D) \times \Cbb$ such that
$\la \psi, 1 \ra=b$,
\beq\label{psig}
\Scal_{\p D}[\psi] + a = f,
\eeq
and
\beq\label{psig2}
\| \psi \|_{-1/2} + |a| \le C (\| f \|_{1/2} + |b|)
\eeq
for some constant $C$. In particular, if we take $b=0$, then $\psi \in H^{-1/2}_0(\p D)$ and \eqnref{psig2} becomes
\beq\label{psig3}
\| \psi \|_{-1/2} + |a| \le C \| f \|_{1/2}.
\eeq
So, we obtain the following lemma.
\begin{lemma}\label{lem:Riesz}
For any $f \in H^{1/2}(\p D)$ there is a unique $\psi \in H^{-1/2}_0(\p D)$ such that
\beq\label{Riesz}
\la \Gvf,f \ra = \la \Gvf, \psi \ra_{\Hcal}
\eeq
for all $\Gvf \in H^{-1/2}_0(\p D)$. Moreover, there is a constant $C$ independent of $f$ such that
\beq\label{Riesz2}
\| \psi \|_{-1/2} \le C \| f \|_{1/2}.
\eeq
\end{lemma}
\pf
For $f \in H^{1/2}(\p D)$ let $\psi \in H^{-1/2}_0(\p D)$ and $a$ be such that \eqnref{psig} and \eqnref{psig3} hold.
Then we have
$$
\la \Gvf, f \ra = \la \Gvf, \Scal_{\p D}[\psi] + a \ra = -\la \Gvf, \psi \ra_{\Hcal}.
$$
Replace $\psi$ by $-\psi$ to have \eqnref{Riesz}.
Uniqueness of $\psi$ is obvious and the proof is complete.
\qed

As a consequence we obtain the following theorem.
\begin{theorem}\label{lemma:norm}
Let $D$ be a bounded Lipschitz domain. There are constants $C_1$ and $C_2$ (which may depend on $D$) such that
\beq\label{Hcal-1/23}
C_1 \| \Gvf \|_{H^{-1/2}(\p D)} \le  \| \Gvf \|_\Hcal  \le C_2\| \Gvf \|_{H^{-1/2}(\p D)}
\eeq
for all $\Gvf  \in H^{-1/2}_0(\p D)$.
\end{theorem}
\pf
Since $\Scal_{\p D}$ is bounded from $H^{-1/2}(\p D)$ into $H^{1/2}(\p D)$, we have for $\Gvf \in \Hcal$
$$
|\la \Gvf,\Gvf \ra_\Hcal | \le \| \Scal_{\p D}[\Gvf] \|_{1/2} \| \Gvf \|_{-1/2} \le C \| \Gvf \|_{-1/2}^2.
$$
 So, we have
\beq\label{Hcal-1/2}
\| \Gvf \|_\Hcal \le C \| \Gvf \|_{-1/2}
\eeq
for some constant $C$.

To prove the opposite inequality, for $f \in H^{1/2}(\p D)$ choose $\psi \in H^{-1/2}_0(\p D)$ satisfying \eqnref{Riesz} and \eqnref{Riesz2}. Then we have, by the Cauchy-Schwarz inequality and \eqnref{Riesz2},
\beq \label{c_eq}
|\la \Gvf,f \ra| = |\la \Gvf, \psi \ra_{\Hcal}| \le \| \Gvf \|_{\Hcal} \| \psi \|_{\Hcal} \le C \| \Gvf \|_{\Hcal} \| \psi \|_{-1/2} \le C
\| \Gvf \|_{\Hcal} \| f \|_{1/2}.
\eeq
Since the above inequality holds for all $f \in H^{1/2}(\p D)$, we obtain
\beq\label{Hcal-1/22}
\| \Gvf \|_{-1/2} \le C \| \Gvf \|_\Hcal .
\eeq
This completes the proof.
\qed

\medskip

Since $\Kcal_{\p D}^*$ is self-adjoint on $\Hcal$, its spectrum $\Gs(\Kcal_{\p D}^*)$ is real, consists of point and continuous spectra, and is a closed set. Moreover, by the spectral resolution theorem (see \cite{Yo80Book}) there is a family of projection operators $\Ecal(t)$ on $\Hcal$ (called a resolution of identity) such that
\beq\label{specresol1}
\Kcal^*_{\p D} = \int_{t \in \Gs(\Kcal_{\p D}^*)} t \, d\Ecal(t).
\eeq
Let $b_{\p D}$ be the spectral bound of $\Kcal^*_{\p D}$, namely,
\beq\label{specbound2}
b_{\p D}:= \sup \{ |\Gl|: \Gl \in \Gs(\Kcal_{\p D}^*) \}.
\eeq
It is well-known that for any Lipschitz domain $b_{\p D} \le \frac{1}{2}$.
In fact, it is proved in \cite{Ke29} that
\beq\label{specbound}
b_{\p D} = \frac{1}{2} \sup_{\Gvf \in \Hcal} \frac{\left| \int_{\Rbb^d \setminus D} |\nabla \Scal_{\p D}[\Gvf]|^2 dx -
\int_{D} |\nabla \Scal_{\p D}[\Gvf]|^2 dx \right|}{\int_{\Rbb^d} |\nabla \Scal_{\p D}[\Gvf]|^2 dx}.
\eeq

We obtain the following lemma.
\begin{theorem}\label{lem:specbound}
Let $D$ be a bounded Lipschitz domain. Then,
\beq\label{specesti}
b_{\p D} < \frac{1}{2}.
\eeq
\end{theorem}
\pf It is proved in \cite{Ve84} that $\Gl I - \Kcal^*_{\p D}$ is invertible on $L^2(\p D)$ if $\Gl \in \mathbb{R}\setminus (-1/2, 1/2]$, and invertible on $L^2_0(\p D)$ if $\Gl \in \mathbb{R}\setminus (-1/2, 1/2)$. Using this result, it is proved in \cite{ChLe08} that $\Gl I -\Kcal^*_{\p D}$ is invertible on $H^{-1/2}(\p D)$ if $\Gl \notin (-1/2, 1/2]$. But, the proof there shows that $\Gl I - \Kcal^*_{\p D}$ is invertible on $H^{-1/2}_0(\p D)$ if $\Gl \notin (-1/2, 1/2)$. We then infer from \eqnref{Hcal-1/23} that the spectrum $\Gs(\Kcal_{\p D}^*)$ of $\Kcal^*_{\p D}$ on $\Hcal$ lies in $(-1/2, 1/2)$. Since $\Gs(\Kcal_{\p D}^*)$ is a closed set, we have \eqnref{specesti}. \qed

Before completing this section let us make a few remarks on the spectrum of the NP operator.
If $\p D$ is $\mathcal{C}^{1, \Ga}$ for some $\Ga >0$, then $\Kcal^*_{\p D}$ is compact and has only a point spectrum accumulating to $0$. The point spectrum of $\Kcal^*_{\p D}$ is completely known when $D$ is a disk, a ball, or an ellipse (see \cite{Kang}). For example, if $D$ is an ellipse of the long axis $a$ and short axis $b$, then eigenvalues of $\Kcal^*_{\p D}$ are
\beq
\pm\frac{1}{2}\left(\frac{a-b}{a+b}\right)^n, \quad n=1,2,\ldots.
\eeq
See \cite{AmKaLe07} \cite{KhPuSh07}. So, if the eccentricity of the ellipse becomes large, the spectral bound approaches to $1/2$.
Other than these examples, the complete spectrum of the NP operator on two discs is obtained in relation to the analysis of gradient concentration \cite{BT} \cite{LiYu14}. On the other hand, a bound for the essential spectrum has been obtained \cite{PePu12}.

\section{Transmission problems in the free space}\label{sec:trans}

The solution $u_k$ to \eqnref{trans} can be represented as
\beq\label{solrep2}
u_k(x) = h(x) + \Scal_{\p D}[\Gvf_k](x), \quad x \in \Rbb^d,
\eeq
where $\Gvf_k \in \Hcal(\p D)$ is the solution to
\beq\label{solinttrans}
\left( \Gl I - \Kcal_{\p D}^* \right)[\Gvf_k] = \p_\nu h |_{\p D} \quad\mbox{on } \p D,
\eeq
where $\Gl=\Gl(k)$ is defined by \eqnref{bilinear}. In fact, \eqnref{solinttrans} is a consequence of the jump relation \eqnref{singlejump} and transmission conditions on $\p D$:
\beq\label{transcon}
u|_+=u|_-, \quad \p_\nu u|_+ = k \p_\nu u|_- \quad\mbox{on } \p D.
\eeq
Note that the righthand side of \eqnref{solinttrans} does not depend on $k$. So, in this case, if $\Gl(k) \notin \Gs(\Kcal_{\p D}^*)$, then
\eqnref{solinttrans} is solvable in $\Hcal$, and hence \eqnref{trans} has a solution. This is already interesting since a negative $k$ can satisfy the condition $\Gl(k) \notin \Gs(\Kcal_{\p D}^*)$. For example, any $k$ with
\beq
k < -\frac{1+2b_{\p D}}{1-2b_{\p D}} \quad\mbox{or}\quad -\frac{1-2b_{\p D}}{1+2b_{\p D}} < k < 0
\eeq
satisfies the condition.

We obtain the following solvability result.

\begin{theorem}\label{thm:solva}
If $\Gl(k) \notin \Gs(\Kcal_{\p D}^*)$, then for any harmonic function $h$ in $\Rbb^d$ there is a unique solution $u_k$ to \eqnref{trans} satisfying
\beq\label{existence}
\| \nabla(u_k - h) \|_{L^2(\Rbb^d)} \le C_k \left\| \p_\nu h \right\|_{-1/2(\p D)}
\eeq
for some $C_k$ depending on $k$.
\end{theorem}

\pf
Note that existence of a solution is already proved. To show \eqnref{existence}, we note that
\begin{align*}
\| \nabla(u_k - h) \|_{L^2(\Rbb^d)}^2 &= \int_{D} |\nabla \Scal_{\p D}[\Gvf_k](x)|^2 dx + \int_{\Rbb^d \setminus D} |\nabla\Scal_{\p D}[\Gvf_k](x)|^2 dx \\
&= \int_{\p D} \p_\nu \Scal_{\p D}[\Gvf_k]|_- \Scal_{\p D}[\Gvf_k] \, d\Gs - \int_{\p D} \p_\nu \Scal_{\p D}[\Gvf_k]|_+ \Scal_{\p D}[\Gvf_k] \, d\Gs \\
&= -\int_{\p D} \Gvf_k \Scal_{\p D}[\Gvf_k] \, d\Gs = \| \Gvf_k \|_{\Hcal}^2.
\end{align*}
So, \eqnref{existence} follows from \eqnref{solinttrans}.

To show uniqueness of the solution, assume that $u_k^1$ and $u_k^2$ are solutions to \eqnref{trans}. Let $v=u_k^1-u_k^2$. Then $v$ is the solution to \eqnref{trans} with $h=0$. So we have
\begin{align*}
0 &= \int_{\Rbb^d} (\chi(\Rbb^d \setminus D) + k \chi(D)) |\nabla v(x)|^2 dx \\
& = \int_{\Rbb^d \setminus D} |\nabla v(x)|^2 dx + k' \int_{D} |\nabla v(x)|^2 dx + i k'' \int_{D} |\nabla v(x)|^2 dx,
\end{align*}
where $k=k'+ik''$. So, if $k'>0$, or if $k' \le 0$ and $k'' \neq 0$, then
$$
\int_{D} |\nabla v(x)|^2 dx= \int_{\Rbb^d \setminus D} |\nabla v(x)|^2 dx=0.
$$
So, $v$ is constant. Since $v(x) \to 0$ as $|x| \to \infty$, we conclude that $v=0$.

Uniqueness for the case $k \le 0$ (and $\Gl(k) \notin \Gs(\Kcal_{\p D}^*)$) can be proved as a limiting case of $k+i\Gd$ as $\Gd \to 0$. For that we need estimate \eqnref{existence} with the constant $C$ independent of $\Gd$. So the proof will be presented after establishing uniformity results in the following.
\qed

\medskip

For $\Ge >0$ let
\beq
R_\Ge:= \{ \, k \in \Cbb \, |\, \mbox{dist} \left( \Gl(k), \Gs(\Kcal_{\p D}^*) \right) \ge \Ge \, \}.
\eeq
It follows from the spectral resolution \eqnref{specresol1} that
\beq\label{specresol2}
\Gvf_k = (\Gl I - \Kcal^*_{\p D})^{-1} \left[ \p_\nu h |_{\p D} \right]  = \int_{t\in \Gs(\Kcal_{\p D}^*)} \frac{1}{\Gl-t} d\Ecal(t) \left[ \p_\nu h|_{\p D} \right] \, ,
\eeq
where $\Gl=\Gl(k)$. If $k \in R_\Ge$, then $|\Gl-t| \ge \Ge$ for all $t \in \Gs(\Kcal_{\p D}^*)$, and hence we obtain from \eqnref{specresol2} that
\begin{align}
\| \Gvf_k \|_{\Hcal}^2 &=  \int_{t\in \Gs(\Kcal_{\p D}^*)} \frac{1}{|\lambda-t|^2}d \left \langle \Ecal(t) \left[ \p_\nu h \right] , \p_\nu h \right \rangle_{\Hcal} \nonumber\\
&\le \frac{1}{\Ge^2} \int_{t\in \Gs(\Kcal_{\p D}^*)} d \left \langle \Ecal(t) \left[ \p_\nu h \right] , \p_\nu h \right \rangle_{\Hcal}
= \frac{1}{\Ge^2} \| \p_\nu h \|_{\Hcal}^2 \, . \label{ortho_est}
\end{align}
Thus we obtain the following theorem from \eqnref{solrep2}.

\begin{theorem}\label{uniitrans}
For each $\Ge>0$ there is a constant $C=C_\Ge>0$ such that
\beq\label{unitrans1}
\| \nabla(u_k - h)\|_{L^2(\Rbb^d)} \le C \left\| \p_\nu h \right\|_{-1/2(\p D)}
\eeq
for all $k\in R_\Ge$ and harmonic functions $h$ in $\Rbb^d$.
\end{theorem}

We also have the following theorem on Lipschitz dependency of the solution on $k$.

\begin{theorem}\label{contitrans}
For any $\Ge>0$ there is a constant $C=C_\Ge>0$ such that
\beq\label{contitrans1}
\| \nabla(u_k - u_s) \|_{L^2(\Rbb^d)} \le \frac{C |k-s|}{(1+|k|)(1+|s|)} \left\| \p_\nu h \right\|_{-1/2(\p D)}
\eeq
for all $k, s \in R_\Ge$ and harmonic functions $h$ in $\Rbb^d$.
\end{theorem}

\pf
We use the spectral resolution \eqnref{specresol2} to have
\begin{align*}
\Gvf_k - \Gvf_s &= \int_{\Gs(\Kcal_{\p D}^*)} \frac{1}{\Gl(k)-t} d\Ecal(t) \left[ \p_\nu h |_{\p D} \right]
- \int_{\Gs(\Kcal_{\p D}^*)} \frac{1}{\Gl(s)-t} d\Ecal(t) \left[ \p_\nu h |_{\p D} \right] \nonumber \\
& = \int_{\Gs(\Kcal_{\p D}^*)} \frac{\Gl(s)-\Gl(k)}{(\Gl(k)-t)(\Gl(s)-t)}d\Ecal(t) \left[ \p_\nu h|_{\p D} \right].
\end{align*}
Note that
$$
\frac{\Gl(s)-\Gl(k)}{(\Gl(k)-t)(\Gl(s)-t)} = \frac{4(k-s)}{[(1-2t)s+1+2t] [(1-2t)k+1+2t]}.
$$
So, if $k, s \in R_\Ge$, we have
\beq
\Big|\frac{\Gl(s)-\Gl(k)}{(\Gl(k)-t)(\Gl(s)-t)}\Big| \le C\frac{|k-s|}{(1+|k|)(1+|s|)}
\eeq
for all $t\in \Gs(\Kcal_{\p D}^*)$. So, we have
$$
\| \Gvf_k - \Gvf_s \|_{-1/2(\p D)} \le   \frac{C|k-s|}{(1+|k|)(1+|s|)} \left\| \p_\nu h \right\|_{-1/2(\p D)} .
$$
Since
$$
u_k(x) - u_s(x) = \Scal_{\p D}[\Gvf_k - \Gvf_s](x), \quad x \in \GO,
$$
we obtain \eqnref{contitrans1}. \qed

\medskip

\noindent{\sl Continuation of proof of Theorem \ref{thm:solva}}.
Suppose $k \le 0$ and $\Gl(k) \notin \Gs(\Kcal_{\p D}^*)$, and let $v=u_k^1-u_k^2$ as before. Choose $\Ge>0$ so that $k \in R_\Ge$. For $\Gd>0$ let $\psi_\Gd$ be the solution to
$$
\left( \frac{k+1+i\Gd}{2(k-1+i\Gd)} I - \Kcal_{\p D}^* \right)[\psi_\Gd] = \frac{i\Gd}{k-1+i\Gd} \p_\nu v|_- \quad\mbox{on } \p D.
$$
Observe that
$$
\frac{k+1+i\Gd}{2(k-1+i\Gd)} = \Gl(k+i\Gd)
$$
and
$$
\left| \frac{k+1+i\Gd}{2(k-1+i\Gd)} - t \right| \ge \Ge \quad\mbox{for all } t \in \Gs(\Kcal_{\p D}^*).
$$
In other words, $k+i\Gd \in R_\Ge$. Then we see from \eqnref{ortho_est} that
\beq\label{psiGd}
\| \psi_\Gd \|_{\Hcal} \le C\Gd \left\| \p_\nu v|_- \right\|_{\Hcal}.
\eeq
Let
\beq
w(x)= v(x) + \Scal_{\p D}[\psi_\Gd](x), \quad x \in \Rbb^d.
\eeq
Then we can see that $w$ satisfies
$$
w|_+=w|_-, \quad \p_\nu w|_+ = (k+i\Gd) \p_\nu w|_- \quad\mbox{on } \p D.
$$
So $w$ is the solution to
$$
 \left \{
 \begin{array}{l}
 \ds \nabla  \cdot ( \chi(\Rbb^d \setminus D) + (k+i\Gd) \chi(D))   \nabla w =0 \quad \mbox{in } \Rbb^d,  \\
 \nm w(x) = O(|x|^{1-d}) \quad\mbox{as } |x| \to \infty.
 \end{array}
 \right.
$$
Since the solution for the complex coefficient $k+i\Gd$ is unique, we obtain from \eqnref{unitrans1} that
$$
\| \nabla w \|_{L^2(\Rbb^d)} =0.
$$
Thus we obtain from \eqnref{psiGd} that
$$
\| \nabla v \|_{L^2(\Rbb^d)} = \| \nabla \Scal_{\p D}[\psi_\Gd] \|_{L^2(\Rbb^d)} \le C \left\| \psi_\Gd \right\|_{\Hcal} \le C\Gd \left\| \p_\nu v|_- \right\|_{\Hcal}.
$$
Letting $\Gd \to 0$, we infer that $v$ is constant in $\Rbb^d$. Since $v(x) \to 0$ as $|x| \to \infty$, we conclude that $v=0$. This completes the proof. \qed

\section{Boundary value problem}\label{sec:bvp}

In this section we consider the boundary value problem \eqnref{NBP}. We first recall a representation formula for a solution $u$ to \eqnref{NBP}.
Suppose that \eqnref{NBP} admits a solution $u \in H^1(\GO)$. Define
the harmonic function $h$ in $\GO$ by
\beq\label{hkx}
h(x)= -\Scal_{\p\GO}[g](x) + \Dcal_{\p\GO}[u|_{\p\GO}](x), \quad x \in \GO,
\eeq
where the double layer potential $\Dcal_{\p\GO}[\Gvf]$ is defined by
\beq
\Dcal_{\p \GO} [\Gvf] (x) := \frac{1}{\Go_d} \int_{\p \GO}
\frac{(y-x)\cdot\nu_y }{|x-y|^d} \Gvf(y)\,d\Gs(y)\;, \quad x\in \Rbb^d\setminus \p \GO.
\eeq
It is proved in \cite{KaSe96} \cite{KaSe00} that $u$ is represented as
\beq\label{solrep}
u(x) = h(x) + \Scal_{\p D}[\Gvf](x), \quad x \in \GO,
\eeq
for some $\Gvf \in \Hcal$. In fact, $\Gvf$ is given by
\beq\label{Gvfform}
\Gvf = \p_\Gv u|_+ - \p_\Gv u|_-
\eeq
and satisfies
\beq\label{solint}
\left( \Gl I - \Kcal_{\p D}^* \right)[\Gvf] = \p_\nu h |_{\p D} \quad\mbox{on } \p D.
\eeq
where $\Gl=\Gl(k)$ is defined by \eqnref{bilinear}.

There is yet another representation for a solution $u$. For $y \in \GO$ let $N_y(x)=N(x,y)$ be the Neumann function on $\GO$, which is the solution to
\beq\label{Neumann}
\begin{cases}
\GD N_y = - \Gd_y \quad & \mbox{in } \GO,\\
\ds \p_\nu N_y = - |\p \GO|^{-1} & \mbox{on } \p \GO,\\
\ds \int_{\p \GO} N_y(x) d \Gs = 0  ,
\end{cases}
\eeq
where $\Gd_y$ is the Dirac mass at $y$ and $|\p \GO|$ denotes the area (or the length) of $\p\GO$. Using $N(x,y)$ we define
\beq
\Ncal_{\p D} [\Gvf](x) = \int_{\p D} N(x,y) \Gvf(y) \, d\Gs(y), \quad x \in \GO.
\eeq
Let  $U$ be the solution in absence of an inclusion, that is, the solution to
\beq\label{PDE}
 \left \{
 \begin{array}{l}
 \ds \GD U =0 \quad \mbox{in } \GO,  \\
 \nm \ds \p_\nu U |_{\p\GO} =g,  \\
 \nm \ds \int_{\p\GO} U \, d\Gs=0.
 \end{array}
 \right .
\eeq
Then it is proved in \cite{AmKa03} (see also \cite{AmKa07Book2}) that with the same $\Gvf$ in \eqnref{Gvfform} the following holds:
\beq\label{Urep}
u(x)=U(x) - \Ncal_{\p D} [\Gvf](x), \quad x \in \p\GO.
\eeq

We prove the following representation theorem.
\begin{theorem}\label{thm:Urep}
Let $u \in H^1(\GO)$ be a solution to \eqnref{NBP}. Then with $\Gvf$ in \eqnref{Gvfform} it holds that
\beq\label{Urep2}
u(x)=U(x) - \Ncal_{\p D} [\Gvf](x), \quad x \in \GO.
\eeq
\end{theorem}

\pf
Let
$$
v(x):= U(x) - \Ncal_{\p D} [\Gvf](x), \quad x \in \GO.
$$
Note that
$$
\p_\Gv \Ncal_{\p D} [\Gvf] = - |\p \GO|^{-1} \int_{\p D} \Gvf =0 \quad \mbox{on } \p\GO.
$$
So, we have $\p_\Gv u= \p_\Gv v$ on $\p\GO$. So by \eqnref{Urep} and unique continuation of harmonic functions, we have
$$
u(x)=v(x), \quad x \in \GO \setminus \overline{D}.
$$

Define $R(x,y)$ by
$$
R(x,y)= N(x,y)+ \GG(x-y),
$$
and let
\beq
\Rcal_{\p D} [\Gvf](x) = \int_{\p D} R(x,y) \Gvf(y) \, d\Gs(y), \quad x \in \GO,
\eeq
so that
\beq\label{NSR}
\Ncal_{\p D}= -\Scal_{\p D} + \Rcal_{\p D}.
\eeq
Since $R(x,y)$ is smooth for $x,y \in \GO$, we have $\p_\Gv \Rcal_{\p D} [\Gvf]|_+ = \p_\Gv \Rcal_{\p D} [\Gvf]|_-$ on $\p D$. Therefore we have on $\p D$
\begin{align*}
\p_\Gv v|_+ - \p_\Gv v|_- &= - \p_\Gv \Ncal_{\p D} [\Gvf]|_+ + \p_\Gv \Ncal_{\p D} [\Gvf]|_- \\
& = \p_\Gv \Scal_{\p D} [\Gvf]|_+ - \p_\Gv \Scal_{\p D} [\Gvf]|_- = \p_\Gv u|_+ - \p_\Gv u|_-.
\end{align*}
So, we have $\p_\Gv v|_-= \p_\Gv u|_-$ on $\p D$. We then infer from unique continuation that $u=v$ in $D$. This completes the proof. \qed

There is a significant difference between representations \eqnref{solrep} and \eqnref{Urep}: $U$ in \eqnref{Urep} is independent of $k$, but $h$ in \eqnref{solrep} depends on $k$ through $u_k|_{\p\GO}$, the Dirichlet data of the solution.

\subsection{Existence and uniqueness of the solution}\label{sec:well}

The purpose of this subsection is to prove the following theorem.

\begin{theorem}\label{thm:well}
There is a compact interval on the negative real axis, say $[a,b]$ with $-\infty < a < b <0$, such that for any $k \in \Cbb \setminus [a,b]$ there is a unique solution $u_k \in H^1(\GO)$ to \eqnref{NBP} such that
\beq\label{ukest}
\| u_k \|_{H^1(\GO)} \le C_k \| g \|_{-1/2(\p\GO)}
\eeq
for some $C_k$ independent of $g$ (which may depend on $k$).
\end{theorem}

\pf
Let $U$ be the solution to \eqnref{PDE}. According to Theorem \ref{thm:Urep}, unique existence of the solution $u$ to \eqnref{NBP} amounts to that of $\Gvf$ in \eqnref{Urep2}, which in turn amounts to unique solvability of the integral equation
$$
k \p_\nu\Ncal_{\p D}[\Gvf]|_- - \p_\nu \Ncal_{\p D}[\Gvf]|_+ = (k-1) \p_\nu U \quad\mbox{on } \p D.
$$
Using \eqnref{singlejump} and \eqnref{NSR} this equation can be written as
$$
k (1/2I - \Kcal_{\p D}^*)[\Gvf] + k \p_\nu\Rcal_{\p D}[\Gvf] + (1/2I + \Kcal_{\p D}^*)[\Gvf] - \p_\nu \Rcal_{\p D}[\Gvf] = (k-1) \p_\nu U,
$$
or equivalently
\beq\label{KReqn}
(\Gl(k) I - \Kcal_{\p D}^*)[\Gvf] + \p_\nu\Rcal_{\p D}[\Gvf] = \p_\nu U.
\eeq\label{KReqn2}
It is convenient to write the above equation as
\beq
(\Gl(k) I - \widetilde{\Kcal}_{\p D}^*)[\Gvf] = \p_\nu U
\eeq
by putting
$$
\widetilde{\Kcal}_{\p D}^* [\Gvf] = \Kcal_{\p D}^*[\Gvf] - \p_\nu\Rcal_{\p D}[\Gvf].
$$

Suppose that $\Gl(k) \in \Cbb \setminus [-b_{\p D}, b_{\p D}]$ where $b_{\p D}$ be the spectral bound of $\Kcal^*_{\p D}$. Since $\Gl(k) I - \Kcal_{\p D}^*$ is invertible on $\Hcal(\p D)$ and $\Gvf \mapsto \p_\nu \Rcal_{\p D}[\Gvf]$ is a compact operator on $\Hcal(\p D)$, unique solvability of \eqnref{KReqn} follows from injectivity. To prove injectivity, suppose that
$$
(\Gl(k) I - \Kcal_{\p D}^*)[\Gvf] + \p_\nu\Rcal_{\p D}[\Gvf] = 0 \quad \mbox{on } \p D.
$$
Let $u(x):= - \Ncal_{\p D}[\Gvf](x)$, $x \in \GO$. Then, $u$ is a solution to \eqnref{NBP} with $g=0$. If $k=0$, then $\p_\nu u|_+=0$ on $\p D$, and hence $u=0$ in $\GO \setminus D$. It then follows that $\Ncal_{\p D}[\Gvf]=0$ in $D$. So we have
$$
\Gvf = \p_\Gv \Ncal_{\p D}[\Gvf] |_- - \p_\Gv \Ncal_{\p D}[\Gvf] |_+ =0.
$$
If $k=\infty$, then $u=\mbox{const}$ on $\p D$. So one can see similarly that $\Gvf=0$.
In general, it can be seen that
$$
k' \int_{D} |\nabla u(x)|^2 dx + \int_{\GO \setminus D} |\nabla u(x)|^2 dx + ik'' \int_{D} |\nabla u(x)|^2 dx =0.
$$
So, if $k' \ge 0$ ($k \neq 0$) or $k'<0$ and $k'' \neq 0$, then $u=0$, and hence $\Gvf =0$.

So far, we have shown that \eqnref{KReqn2} has a unique solution if $k \in \Cbb \setminus (-\infty, 0)$.
Suppose that $|k| >R$ for some $R$ to be determined. Note that $\Gl(\infty)=1/2$. We then have
$$
|\Gl(k)-1/2|= \left| \frac{k+1}{2(k-1)} - \frac{1}{2} \right| = \frac{1}{|k-1|} \le \frac{1}{R-1}.
$$
If $R$ is sufficiently large, it follows that
\begin{align*}
\Gl(k) I - \widetilde{\Kcal}_{\p D}^* & = \left(\Gl(k)- 1/2 \right)I + \left( 1/2I - \widetilde{\Kcal}_{\p D}^* \right) \\
&= \left( 1/2 I - \widetilde{\Kcal}_{\p D}^* \right) \left[ I+ \left(\Gl(k)-1/2 \right) \left( 1/2I - \widetilde{\Kcal}_{\p D}^* \right)^{-1} \right],
\end{align*}
and hence there is $C$ such that
\beq\label{unifest}
\| (\Gl(k) I - \widetilde{\Kcal}_{\p D}^*)^{-1} [\Gvf] \|_{1/2(\p D)} \le C\| \Gvf \|_{1/2(\p D)}
\eeq
for all $|k| > R$. Similarly, one can show that there is $r<1$ such that \eqnref{unifest} holds for all $|k|<r$.

Let $a=-R$ and $b=-r$. We have shown that \eqnref{KReqn2} has a unique solution if $k \in \Cbb \setminus [a,b]$ for some $-\infty < a < b <0$. Let $\Gvf_k$ be the solution, and define
$$
u_k(x)=U(x) - \Ncal_{\p D} [\Gvf_k](x), \quad x \in \GO.
$$
Then, $u_k$ is the solution to \eqnref{NBP} satisfying \eqnref{ukest}. \qed

During the course of proof above we proved the following lemma.
\begin{lemma}\label{lem:ukest2}
There is $C$ (independent of $k$ and $g$) such that
\beq\label{ukest2}
\| u_k \|_{H^1(\GO)} \le C \| g \|_{-1/2(\p\GO)}
\eeq
for all $k$ satisfying $|k| >-a$ or $|k| <-b$.
\end{lemma}

The following example shows that for the boundary value problem uniqueness may fail, unlike the free space problem in the previous section, even if $\Gl(k) \notin \Gs(\Kcal_{\p D}^*)$.

\begin{example}
Let $D$ be the disk centered at $0$ of radius $r_i$. Then $\Kcal_{\p D}^*=0$. So, the free space problem \eqnref{trans} admits a unique solution for all $k \neq -1$ ($k \in \Cbb$). But it is not the case for the boundary value problem. Let $\GO$ be the disk centered at $0$ of radius $r_e$ ($r_i < r_e$). Let $\rho:= r_i/r_e$. For a positive integer $n$ let
$$
k_n:= \frac{\rho^{2n}-1}{\rho^{2n}+1}.
$$
Then the function $v$, defined by
$$
v(x) :=
\begin{cases}
(r^n+ r_e^{2n} r^{-n}) e^{in\Gt}, \quad &\mbox{if } r_i < r=|x| \le r_e, \\
\nm
\ds \left(1+ \rho^{-2n} \right) r^n e^{in\Gt}, \quad &\mbox{if }  r \le r_i,
\end{cases}
$$
is a solution to \eqnref{NBP} with $k=k_n$. Note that $\p_\nu v=0$ on $\p \GO$. It is quite interesting to observe that $\frac{1}{2} \rho^{2n}$ is the eigenvalues of the NP operator associated with two interfaces $\p D$ and $\p\GO$ (see {\rm\cite{ACKLM}}).
\end{example}

\subsection{Uniformity of regularity estimates}\label{sec:uni}

In this section and sections to follow we use the condition \eqnref{kcondition}. This condition is also required for estimation of the Neumann-to-Dirichlet map (see Lemma \ref{Di-N}). For a positive constant $L>0$ let $S_L$ be the infinite sector of $k$ satisfying \eqnref{kcondition} and let $O_L$ be the image of $S_L$ under the transformation \eqnref{bilinear} (the grey region in the right figure in Figure \ref{fig}).
Note that for a given $L>0$ there is $\Ge>0$ such that $S_L \subset R_\Ge$.
In fact, we have
\beq
\mbox{dist}(O_L, [-b_{\p D}, b_{\p D}] ) = 2\left(\frac{1}{4} - b_{\p D}^2\right) \left(\sqrt{1 + L^{-2}} + \sqrt{(2 b_{\p D})^2 + L^{-2}}\right)^{-1}.
\eeq
It is worth mentioning that $\mbox{dist}(O_L, [-b_{\p D}, b_{\p D}] ) \to 0$ if either the spectral bound $b_{\p D}$ tends to $1/2$ or $L \to 0$.

\begin{lemma}\label{Di-N}
For each $L>0$ there is a constant $C=C_L$ such that
\beq\label{ukpgo}
\| u_k \|_{1/2(\p\GO)} \le C \| g \|_{-1/2(\p\GO)}
\eeq
for all $k \in S_L$ and $g \in H^{-1/2}_0(\p\GO)$.
\end{lemma}

\pf Let $\la \ , \ \ra_{\p \GO}$ be the $H^{-1/2}$-$H^{1/2}$ pairing on $\p\GO$. Note that
$$
\la g, u_k \ra_{\p\GO} = \int_{\GO} \Gg_k |\nabla u_k|^2 dx  = \int_{\GO\setminus D} |\nabla u_k|^2 dx + k' \int_{D} |\nabla u_k|^2 dx + i k'' \int_{D} |\nabla u_k|^2 dx .
$$
So we have
\beq \label{pre-ukpgo3}
|\la g, u_k \ra_{\p\GO} |^2   = \Big(\|\nabla u_k\|_{L^2(\GO \setminus D)}^2 + k' \|\nabla u_k\|_{L^2(D)}^2\Big)^2 + k''^2 \|\nabla u_k\|_{L^2(D)}^4.
\eeq

If $k' \geq 0$, then
\beq\label{kposi}
\|\nabla u_k\|_{L^2(\GO \setminus  D)}^2 \le |\la g, u_k \ra_{\p\GO} |  \le \|g\|_{-1/2(\p \GO)} \|u_k\|_{1/2(\p \GO)}.
\eeq
If $k' < 0$, we rewrite \eqnref{pre-ukpgo3} as
\beq \label{pre-ukpgo4}
|\la g, u_k \ra_{\p\GO} |^2   = |k|^2 \Big( \|\nabla u_k\|_{L^2(D)}^2 + \frac{k'}{|k|^2} \|\nabla u_k\|_{L^2(\GO \setminus D)}^2 \Big)^2 + \frac{k''^2}{|k|^2} \|\nabla u_k\|_{L^2(\GO \setminus D)}^4.
\eeq
Since $k' <0$, the first term on the righthand side can be $0$. On the other hand, $k''^2/|k|^2> C$ for some constant $C>0$ if and only if
$k''^2 \geq  L^2 k'^2$ for some $L$. So, we obtain that if $k \in S_L$, then
\beq\label{knega}
\|\nabla u_k\|^2_{L^2(\GO \setminus D)}  \leq C \|g\|_{-1/2(\p \GO)} \|u_k\|_{1/2(\p \GO)}
\eeq
for some constant $C$ independent of $k$.

Choose a smooth subdomain $\GO_0$ of $\GO$ containing $\overline{D}$. For $\eta \in H^{-1/2}_0(\p \GO)$ let $w$ be the solution to
 \begin{equation*}
 \left \{
 \begin{array}{l}
 \ds \GD w =0 \quad \mbox{in } \GO \setminus \overline{\GO_0},  \\
 \nm \ds \p_\nu w |_{\p\GO} = {\eta}, \quad  \p_\nu w |_{\p \GO_0} = 0,\\
 \nm \ds\int_{\GO \setminus {\GO_0}} w(x) dx  =0.
 \end{array}
 \right .
 \end{equation*}
Then there is a constant $C$ such that
$$
\|\nabla w \|_{L^2 (\GO \setminus \GO_0)} \leq C \|\eta\|_{-1/2(\p \GO)}.
$$
Since
\begin{equation*}
\la \eta, u_k \ra_{\p\GO}  =  \int_{\GO \setminus \GO_0} \nabla w \cdot \nabla \overline{u_k} \, dx,
\end{equation*}
we have
\begin{equation*}
|\la \eta, u_k \ra_{\p\GO} |  \le C \|\nabla u_k\|_{L^2 (\GO \setminus \GO_0)} \|\eta\|_{-1/2(\p \GO)},
\end{equation*}
and hence
\beq\label{pre-ukpgo}
\| u_k \|_{1/2(\p\GO)} \le C \|\nabla u_k\|_{L^2 (\GO \setminus \GO_0)}.
\eeq
Now, \eqnref{ukpgo} follows from \eqnref{knega} and \eqnref{pre-ukpgo}.
\qed

As a consequence we obtain the following lemma.
\begin{lemma}\label{thm:uni}
For any $L>0$ there is a constant $C=C_L$ such that
\beq\label{ukh11}
\| u_k \|_{H^1(\GO)} \le C \| g \|_{-1/2(\p\GO)}
\eeq
for all $k \in S_L$ and $g \in H^{-1/2}_0(\p\GO)$.
\end{lemma}

\pf
By the same way as in \eqnref{ortho_est} we obtain
\beq\label{ortho_est2}
\| \Gvf_k \|_{\Hcal} \le C \left\| \p_\nu h_k \right\|_{\Hcal} \, .
\eeq
We then obtain from \eqnref{Hcal-1/23} that
\beq\label{gvfh}
\| \Gvf_k \|_{-1/2(\p D)} \le C \left\|\p_\nu h_k \right\|_{-1/2(\p D)}
\eeq
for some constant $C$ independent of $k \in S_L$. Here (and throughout this paper) the constant $C$ may differ at each occurrence.

It follows from \eqnref{hkx} and \eqnref{ukpgo} that
\begin{align}
\| h_k \|_{H^1(\GO)} & \le \| \Scal_{\p\GO}[g] \|_{H^1(\GO)} + \| \Dcal_{\p\GO}[u_k|_{\p\GO}] \|_{H^1(\GO)} \nonumber \\
& \le C (\| g \|_{-1/2(\p\GO)} + \| u_k \|_{1/2(\p\GO)}) \le C \| g \|_{-1/2(\p\GO)} \label{hkg}
\end{align}
for some $C$ independent of $k$.  As a consequence, we obtain
\beq\label{hest}
\left\| \p_\nu h_k \right\|_{-1/2(\p D)} \le C \| h_k \|_{H^1(\GO)} \le C \| g \|_{-1/2(\p\GO)}.
\eeq
So, we have from \eqnref{gvfh}
\beq\label{hkg2}
\| \Scal_{\p D}[\Gvf_k] \|_{H^1(\GO)} \le C \| \Gvf_k \|_{-1/2(\p D)} \le C \| g \|_{-1/2(\p\GO)}.
\eeq
Then, \eqnref{ukh11} follows from \eqnref{solrep}, \eqnref{hkg} and \eqnref{hkg2}.
\qed

Let $a$ and $b$ be numbers appearing in Theorem \ref{thm:well} and Lemma \ref{lem:ukest2}. For $\Ge >0$, define
\beq
Q_\Ge: = \{ k \in \Cbb : \text{dist} (k, [a,b])> \Ge \}.
\eeq
For a given $\Ge$ there is $L$ such that
$$
Q_\Ge \subset \{ |k| > -a \} \cup \{ |k| < -b \} \cup S_L.
$$
So we obtain the following theorem from Lemma \ref{lem:ukest2} and Lemma \ref{thm:uni}.
\begin{theorem}\label{thm:uni2}
For any $\Ge>0$ there is a constant $C=C_\Ge$ such that
\beq\label{ukh110}
\| u_k \|_{H^1(\GO)} \le C \| g \|_{-1/2(\p\GO)}
\eeq
for all $k \in Q_\Ge$ and $g \in H^{-1/2}_0(\p\GO)$.
\end{theorem}

\subsection{Lipschitz continuity estimate with respect to the conductivity}\label{sec:con}

In this section we investigate the continuous dependency of $u_k$ on the conductivity parameter $k$. There has been some work on this problem. If $k$ is real and approaches to $\infty$ (or $0$), then it was proved that $u_k$ converges to $u_\infty$ (or $u_0$) in $H^1$-norm (\cite{FrVo89}) and in $L^\infty$-norm (\cite{KaSe99}).

The solution $u_k$ to \eqnref{NBP} depends on $k$ analytically. In fact, one can see easily that $\bar\p_k u_k= \pd{u_k}{\bar k}$ is the solution to
\beq\label{NBP2}
 \left \{
 \begin{array}{l}
 \ds \nabla  \cdot ( \Gg_k  \nabla u) =0 \quad \mbox{in } \GO,  \\
 \nm \ds \p_\nu u |_{\p\GO} =0, \\
 \nm \ds \int_{\p\GO} u \, d\Gs=0,
 \end{array}
 \right.
 \eeq
and hence $\bar\p_k u_k=0$ for $k \in \Cbb \setminus [a,b]$. The following theorem quantifies Lipschitz dependency of the solution on $k$.

\begin{theorem}\label{conti2}
For each $\Ge>0$ there is a constant $C=C_\Ge$ such that
\beq\label{ukus1}
\| u_k - u_s \|_{H^1(\GO)} \le C \frac{|k-s|}{(1+{|k|})(1+{|s|})} \| g \|_{-1/2}
\eeq
for all $k, s \in Q_\Ge$ and $g \in H^{-1/2}_0(\p\GO)$.
\end{theorem}

Before proving Theorem \ref{conti2}, we make a brief remark on the expression of \eqnref{ukus1}. If $k$ and $s$ are bounded, then it is a Lipschitz continuity estimate. But if $k$ and $s$ are large, then it means more than Lipschitz continuity. It shows, for example, that $u_k-u_s$ can be arbitrarily small in $H^1$-norm if  $k$ and $s$ are large.

\medskip

\noindent{\sl Proof of Theorem \ref{conti2}}.
Suppose $k, s \in Q_\Ge$.
Let us define $\tilde u_s$ by
\beq\label{solrepf2}
\tilde u_s(x) := h_k(x) + \Scal_{\p D}[\tilde\Gvf_s](x), \quad x \in \GO,
\eeq
where $\tilde\Gvf_s \in H^{-1/2}_0(\p D)$ is the solution to
\beq\label{solint2}
\left( \Gl(s) I - \Kcal_{\p D}^* \right)[\tilde\Gvf_s] = \p_\nu h_k |_{\p D} \quad\mbox{on } \p D.
\eeq
It is worth mentioning that $\tilde u_s$ is a solution to the equation $\nabla \cdot (\Gg_s \nabla u)=0$ in $\GO$.

We now compare $u_k$ with $\tilde u_s$, and then $\tilde u_s$ with $u_s$.
We use the spectral resolution \eqnref{specresol2} to have
$$
\Gvf_k - \tilde \Gvf_s = \int_{\Gs(\Kcal_{\p D}^*)} \frac{1}{\Gl(k)-t} d\Ecal(t) \left[ \p_\nu h_k |_{\p D} \right]
- \int_{\Gs(\Kcal_{\p D}^*)} \frac{1}{\Gl(s)-t} d\Ecal(t) \left[ \p_\nu h_k |_{\p D} \right] .
$$
So using the same argument as in the proof of Theorem \ref{contitrans} and \eqnref{hest}, we see that
\beq \label{serf56}
\| \Gvf_k - \tilde\Gvf_s \|_{-1/2(\p D)} \le   \frac{C|k-s|}{(1+|k|)(1+|s|)} \left\| \p_\nu h_k \right\|_{-1/2(\p D)} \le   \frac{C|k-s|}{(1+|k|)(1+|s|)} \left\| g \right\|_{-1/2(\p \GO)}.
\eeq
Since
$$
u_k(x) - \tilde u_s(x) = \Scal_{\p D}[\Gvf_k - \tilde\Gvf_s](x), \quad x \in \GO,
$$
we have
\beq \label{k_s_est}
\| u_k - \tilde u_s \|_{H^1(\GO)} \le  \frac{C|k-s|}{(1+|k|)(1+|s|)} \left\| g \right\|_{-1/2(\p \GO)}. \eeq

On the other hand, we obtain from Theorem \ref{thm:uni}
\beq
\|  \tilde u_s -u_s\|_{H^1(\GO)} \le C\left\| \p_\nu \tilde u_s -  \p_\nu u_k \right\|_{-1/2(\p \GO)},
\eeq
and
\begin{align*}
\left\| \p_\nu \tilde u_s -  \p_\nu u_k \right\|_{-1/2(\p \GO)} = \left\| \p_\nu \Scal_{\p D} [\tilde \Gvf_s -\Gvf_k] \, \right\|_{-1/2(\p \GO)} \le C \left\| \tilde \Gvf_s -\Gvf_k \, \right\|_{-1/2(\p D)} ,
\end{align*}
where the last inequality holds since there is a distance between $\p D$ and $\p\GO$. We then infer from \eqnref{serf56} that
\beq \label{ts_s_est}
\|  \tilde u_s -u_s\|_{H^1(\GO)} \le  \frac{C|k-s|}{(1+|k|)(1+|s|)} \left\| g \right\|_{-1/2(\p \GO)}.
\eeq
We obtain \eqnref{ukus1} from \eqnref{k_s_est} and \eqnref{ts_s_est}.
\qed

\subsection{Dirichlet problem}

Let us briefly mention on the Dirichlet boundary value problem:
\beq\label{DBP}
 \left \{
 \begin{array}{l}
 \ds \nabla  \cdot ( \Gg_k  \nabla v) =0 \quad \mbox{in } \GO,  \\
 \nm \ds v |_{\p\GO} =f, \\
 \nm \ds \int_{\p\GO} v \, d\Gs=0.
 \end{array}
 \right.
 \eeq
The representation formula \eqnref{solrep} for the solution $v_k$ is still valid if we use $h_k$ with $g$ replaced with $\p_\nu v_k|_{\p\GO}$ and $u_k|_{\p\GO}$ replaced with $f$ in \eqnref{hkx}. Similarly to Lemma \ref{Di-N} one can show that
$$
\| \p_\nu v_k \|_{-1/2(\p\GO)} \le C \| f \|_{1/2(\p\GO)}
$$
for all $k \in Q_\Ge$.
So following the same lines of proofs we can obtain results for the Dirichlet problem similar to Theorem \ref{thm:uni} and Theorem \ref{conti2}.

\section{Uniform validity of the boundary perturbation formula}\label{sec:asymp}

Let  $U$ be the solution to \eqnref{PDE}. Then $u_k-U$ can be regarded as a perturbation due to presence of the inclusion with the conductivity $k$. If $D$ is diametrically small, then the perturbation is small and the asymptotic formula as the diameter tends to zero is known.
We assume that the inclusion $D$ is represented as
\begin{equation}\label{domain D}
D= \Gd B + z
\end{equation}
where $\Gd$ represents the small diameter of $D$, $B$ is a reference domain containing $0$, and $z$ is the location of $D$. Here we also assume that $D$ is away from $\p\GO$, namely, there is $c_0>0$ such that
\beq
\mbox{dist}(D, \p\GO) \ge c_0.
\eeq
This condition is required so that the interaction between $D$ and $\p \GO$ does not appear in the boundary perturbation formula in the following. Then, the following asymptotic expansion of the boundary perturbation holds:
\beq\label{firstorder}
u_k (x) - U(x) = -\Gd^d \nabla U(z) \cdot M \nabla_z N(x,z) + E_k^1(x), \quad x \in \p\GO,
\eeq
where $M=M(k,B)$ is the polarization tensor associated with the conductivity $k$ and the domain $B$ (see \eqnref{GPT} below for the definition of the polarization tensor), $N(x,z)$ is the Neumann function \eqnref{Neumann},
and $E_k^1(x)$ indicates the error of the approximation and satisfies
\beq\label{errorone}
|E_k^1(x)| \le C_k \Gd^{d+1}
\eeq
for some $C_k$ which may depend on $k$. The formula \eqnref{firstorder} was first discovered in \cite{FrVo89} and used effectively for the inverse problem to find the location and/or some geometric properties (especially the equivalent ellipse of the inclusion) of the inclusion using the boundary measurements \cite{AKKL05} \cite{bhv03} \cite{CFMoVo98}.

The formula has been extended in \cite{AmKa03} to include the higher order terms as
\beq\label{higherorder}
u_k (x) - U(x) = -\sum_{n=2}^{d+1} \sum_{|\Ga| + |\Gb|=n} \frac{\Gd^{n+d-2}}{\Ga!\Gb!} \p^\Ga U(z) m_{\Ga\Gb} \p_z^\Gb N(x,z) + E_k(x), \quad x \in \p\GO,
\eeq
where $m_{\Ga \Gb}= m_{\Ga \Gb}(k,B)$ is a series of tensors, called generalized polarization tensors (GPTs), associated with $(k,B)$ (see \eqnref{GPT} below). Here $\Ga$ and $\Gb$ are multi-indices. The error of this approximation satisfies
\beq\label{errortwo}
|E_k(x)| \le C_k \Gd^{2d}.
\eeq
The formula \eqnref{higherorder} has been also used to solve inverse problems for which we refer readers to \cite{AmKa04Book1} \cite{AmKa07Book2} and references therein. It is worth mentioning that the asymptotic expansion method has been applied in various contexts such as multi-static imaging and bio-medical imaging for which we refer to \cite{AGJ} \cite{AK11} and references in. The first order formula \eqnref{firstorder} was also generalized to the case when the inclusion is an arbitrary subset of $\GO$ of low volume fraction \cite{cv03}. As pointed out in the paper, there may not be a higher order formula in such a case.

Since the conductivity of the inclusion can be extreme (close or equal to $0$ or $\infty$) or complex, it is important to clarify dependence of the constant $C_k$ appearing in \eqnref{errorone} and \eqnref{errortwo}. Recently it is proved in \cite{NgVo09} that the approximation formula \eqnref{firstorder} is valid uniformly in $k$ for $k$ real and $0 \le k \le \infty$, namely, $C_k$ in \eqnref{errorone} can be chosen independently of $k$. The purpose of this section is to show that the approximation formula \eqnref{higherorder} is valid uniformly in complex $k \in Q_\Ge$. It is worth mentioning that an asymptotic formula for elasticity similar to \eqnref{higherorder} was obtained in \cite{AKNT02} and its uniform validity for the real Lam\'e parameters was proved in \cite{AmKaKiLe13} (using a variational method).

Let us recall the definition of GPTs associated to the inclusion $B$ with the conductivity $k$. For a given multi-index $\Ga\in \mathbb{N}^d$, let $\Gvf_{k,\Ga}$ be the solution to
$$
\left( \Gl(k) I - \Kcal_{\p D}^* \right)[\Gvf_{k,\Ga}] = \p_\nu x^\Ga \quad\mbox{on } \p D.
$$
Here $x^{\Ga} = x_1^{\Ga_1} \cdots x_d^{\Ga_d}$. We see through the same estimates as in \eqnref{ortho_est} that there is a constant $C$ such that
\beq \label{phi_bdd}
\|\Gvf_{k,\Ga}\|_{H^{-1/2}(\p B)} \leq C
\eeq
for all $k\in Q_\Ge$. The GPT associated with $(k, B)$ is defined by
\beq\label{GPT}
m_{\Ga \Gb} (k,B) = \int_{\p B} x^{\Gb} \Gvf_{k,\Ga}(x) d \Gs
\eeq
for $\Ga, \Gb \in \mathbb{N}^d$. Here the integral is understood as the $H^{1/2}$-$H^{-1/2}$ pairing.
Thus GPTs for the Lipschitz domain $B$ are bounded independently of $k \in Q_\Ge$.

We obtain the following theorem. The proof is nothing but repetition of the arguments to derive \eqnref{higherorder} in \cite{AmKa03} except that we need to keep track of the dependency on $k$ of the constants appearing during the derivation using results of this paper. So we omit the proof.

\begin{theorem}\label{uniform}
Let $E_k$ be the error defined in \eqnref{higherorder}. For each $\Ge$ there is a constant $C=C_\Ge$ such that
\beq\label{eklinfty1}
\| E_k \|_{L^\infty(\p\GO)} \le C\delta^{2d} \| g \|_{-1/2(\p \GO)}.
\eeq
for all $k \in Q_\Ge$ and $g \in H^{-1/2}_0(\p\GO)$.
\end{theorem}

\section{Multiple inclusions}\label{sec:multi}

In this section, we consider the case when there are multiple inclusions $D_1, \ldots, D_M$ whose closures are mutually disjoint. As before, $D_j$'s are simply connected and have Lipschitz boundaries. We assume that the inclusions are at some distance from $\p\GO$, that is,
there is a constant $c_0$ such that
\beq
\mbox{dist}(D_j, \p \GO) \geq c_0.
\eeq
Suppose that $D_j$ has the complex conductivity $k_j$ ($j=1, \ldots, M$) so that the conductivity distribution is given by
\beq
\Gg_\Bk = \chi(\GO \setminus \cup_{j=1}^M  D_j) + \sum_{j=1}^M k_j \chi(D_j).
\eeq
Here $\Bk=(k_1, \ldots, k_M)$. We seek conditions on $\Bk$ which is sufficient for solvability and uniformity of the estimates for the problem
\beq\label{NBPmulti}
 \left \{
 \begin{array}{l}
 \ds \nabla  \cdot ( \Gg_\Bk  \nabla u) =0 \quad \mbox{in } \GO,  \\
 \nm \ds \p_\nu u |_{\p\GO} =g, \\
 \nm \ds \int_{\p\GO} u \, d\Gs=0.
 \end{array}
 \right.
 \eeq

\subsection{The NP operator and symmetrization}

Let us use notation $\Scal_j$, $\Kcal_j^*$ and $\Ncal_j$ for $\Scal_{\p D_j}$, $\Kcal_{\p D_j}^*$ and $\Ncal_{\p D_j}$ to make expressions simpler. Similarly to \eqref{solrep}, the solution $u_{\Bk}$ to  \eqref{NBPmulti} can be uniquely represented as
\beq\label{ukx}
u_{\Bk}(x) = h_{\Bk}(x) + \sum_{j=1}^M \Scal_{j} [\Gvf^\Bk_j] (x), \quad x \in \GO,
\eeq
where the harmonic function $h_{\Bk}$ is given by \eqref{hkx} and $\Gvf^\Bk_j \in H^{-1/2}_0(\p D_j)$, $j = 1, \cdots, M$, is the solution to the system of the integral equations
\begin{equation*}
(\Gl_j I - \Kcal_{j}^*) [\Gvf^\Bk_j] - \sum_{l\neq j} \p_{\nu_j} \Scal_{l} [\Gvf^\Bk_l] = \p_{\nu_j} h_{\Bk} \quad \mbox{on } \p D_j, \quad \Gl_j = \Gl(k_j),
\end{equation*}
where $\p_{\nu_j}$ denotes the outward normal derivative on $\p D_j$. Set
\beq
\Dbb (\GL) = \Dbb(\Gl_1, \ldots, \Gl_M) := \mbox{diag}[\Gl_1, \ldots, \Gl_M] ,
\eeq
and
\begin{equation*}
\Phi_\Bk := [\Gvf^\Bk_1, \ldots, \Gvf^\Bk_M]^T, \quad  \p h_\Bk = \left[\p_{\nu_1} h_\Bk |_{\p D_1}, \ldots, \p_{\nu_M} h_\Bk |_{\p D_M}\right]^T \ (T \mbox{ for transpose}).
\end{equation*}
Then the above system of integral equations can be rewritten in a matrix form
\beq\label{solmatrix}
\left(\Dbb(\GL) - \Kbb^*\right)[\Phi_\Bk] = \p h_\Bk,
\eeq
where
\beq\label{NPmulti}
 \Kbb^* := \begin{bmatrix} \Kcal_{1}^* & \p_{\nu_1} \Scal_{2} & \cdots & \p_{\nu_1} \Scal_{M}\\
\p_{\nu_2} \Scal_{1} & \Kcal_{2}^* & \cdots & \p_{\nu_2} \Scal_{M} \\
\vdots & \vdots & \ddots & \vdots\\
\p_{\nu_M} \Scal_{1} & \p_{\nu_M} \Scal_{2} & \cdots & \Kcal^*_{M}\end{bmatrix}.
\eeq
The operator $\Kbb^*$ is the NP operator corresponding to multiple inclusions $(D_1, \ldots, D_M)$.

The equation \eqnref{solmatrix} holds in the space $\BH := \prod_{j=1}^M H^{-1/2}_0(\p D_j)$. As in the single interface case, $\Kbb^*$ is not self-adjoint on $\BH$, but it can be symmetrized by introducing a new inner product on $\BH$. Define, for $i,j=1, \ldots, M$, $\Scal_{ij}: H^{-1/2}_0(\p D_j) \to H^{1/2}_0(\p D_i)$ by
\beq
\Scal_{ij}[\Gvf](x):= \Scal_j[\Gvf](x), \quad x \in \p D_i \, ,
\eeq
and $\Sbb$ on $\mathbf{H}$ by
\begin{equation*}
 \Sbb :=
  \begin{bmatrix} \Scal_{11} &\Scal_{12} & \cdots & \Scal_{1M} \\
  \Scal_{21} &\Scal_{22} & \cdots & \Scal_{2M} \\
\vdots & \vdots & \ddots & \vdots\\
  \Scal_{M1} &\Scal_{M2} & \cdots & \Scal_{MM}
  \end{bmatrix}.
\end{equation*}
The twisted inner product \eqnref{innerp} is now extended to
\beq\label{innermulti}
\la \Phi, \Psi \ra_\Hcal := - \la \Phi, \Sbb[\Psi] \ra, \quad \Phi, \Psi \in \mathbf{H}.
\eeq
Then one can show that
\begin{itemize}
\item[(i)] $\la \Phi, \Psi \ra_\Hcal$ is an inner product on $\BH$.
\item[(ii)] The symmetrization principle \eqnref{Plemelj} can be extended as
\beq
\Sbb \Kbb^* = \Kbb \Sbb,
\eeq
where $\Kbb$ is the adjoint of $\Kbb^*$ with respect to $L^2$ inner product.
\item[(iii)] The spectrum of $\Kbb^*$ on $\BH$ lies in $(-1/2, 1/2)$.
\end{itemize}
In fact, these facts were proved in \cite{ACKLM} where there are two interfaces. Following the same lines of proofs there one can extend them to the case of multiple inclusions.
Let $\| \, \cdot\,\|_\Hcal$ be the norm induced by the inner product $\la \cdot, \cdot  \ra_\Hcal$. Then we have the following lemma.

\begin{lemma}
There are constants $C_1$ and $C_2$ such that
\beq\label{norm_equiv}
 C_1 \|\Phi\|_{\BH} \leq \|\Phi\|_{\Hcal} \leq C_2  \|\Phi\|_{\BH} \quad\mbox{for all } \Phi \in \BH.
\eeq
\end{lemma}
\pf
Similarly to the operator $A$ in Lemma \ref{lemA} we define the operator
$$
A: \prod_{j=1}^M H^{-1/2}(\p D_j) \times \Cbb^M \to \prod_{j=1}^M H^{1/2}(\p D_j) \times \Cbb^M
$$ by
\beq\label{defAm}
A(\Phi, (a_1, \ldots, a_M)) := \left( \Sbb[\Phi] + (a_1, \ldots, a_M), (\la \Gvf_1, 1 \ra_1, \ldots, \la \Gvf_M, 1 \ra_M)  \right)
\eeq
where $\Phi= [\Gvf_1, \ldots,\Gvf_M]^T$ and $\la \Gvf_j, 1 \ra_j$ denotes $H^{-1/2}$-$H^{1/2}$ duality pairing on $\p D_j$. One can show following the same lines of the proof of Lemma \ref{lemA} that $A$ is invertible. Then \eqnref{norm_equiv} follows by the same argument as in section \ref{sec:spec}. \qed

\subsection{A complete characterization of solvability for two disks}\label{sec:twodisk}

Suppose that multiple inclusions consist of two disks $D_1$ and $D_2$ of the same radius. Suppose that the conductivity of $D_j$ is $k_j$ for $j=1,2$ so that the conductivity distribution is given by
\beq\label{sigma:def}
\Gg=k_1\chi(D_1)+k_2\chi(D_2)+\chi(\Rbb^2\setminus (D_1\cup D_2)).
\eeq
We consider solvability in terms of $k_1$ and $k_2$ of the problem
\beq\label{cond_eqn}
\begin{cases}
\ds\nabla\cdot\Gg\nabla u=0\quad&\mbox{in }\Rbb^2, \\
\ds u(\Bx) - h(\Bx)  =O(|\Bx|^{-1}) \quad&\mbox{as } |\Bx| \to \infty,
\end{cases}
\end{equation}
where $h$ is a given harmonic function in $\Rbb^2$.

For that purpose we use the bipolar coordinates.
The bipolar coordinate system $(\xi,\Gt)\in\Rbb\times (-\pi,\pi]$ is defined as follows (see \cite{Lock}): for a fixed $\Ga >0$
\beq  \label{bipolar}
x=\Ga\frac{ \sinh \xi }{\cosh \xi - \cos \Gt} \quad \mbox{ and } \quad y=\Ga\frac{\sin \Gt}{\cosh \xi - \cos \Gt}.
\eeq
With an appropriately chosen $\Ga$ two disks $D_1$ and $D_2$ are given in terms of bipolar coordinates by
\beq
D_1=\{\xi<-\xi_0\},\quad D_2=\{\xi>\xi_0\}
\eeq
for some $\xi_0>0$. Recall that the normal derivative of a function $u$ on $\p D_j$ are given by
\beq\label{nor_bipolar}
\pd{u}{\nu}\Bigr|_{\xi=c}=-\mbox{sgn}(c)J(c,\Gt)\pd{u}{\xi}\Bigr|_{\xi=c},
\eeq
where
$$
J(\xi,\Gt):=\frac{\cosh\xi-\cos\Gt}{\Ga}.
$$

The bipolar coordinate system admits a general separation of variables solution to the harmonic function $u$ as follows:
\begin{align}
u(\xi,\Gt)&=a_0+b_0\xi+c_0\Gt+\sum_{n=1}^\infty\bigr[ (a_n e^{n\xi}+b_n e^{-n\xi})\cos n\Gt+\bigr (c_n e^{n\xi}+d_n e^{-n\xi})\sin n\Gt\bigr],
\end{align}
where $a_n$, $b_n$, $c_n$ and $d_n$ are constants.
Suppose that the expansion of the harmonic function $h$ is given by
\begin{align}
h(\xi,\Gt)&=f_0+\begin{cases}
\ds\sum_{n\neq 0} f_n e^{-|n|\xi}e^{i n\Gt},& \xi>0,\\
\ds
\sum_{n\neq 0} g_n e^{|n|\xi}e^{i n\Gt},&\xi<0,
\end{cases}
\end{align}
for some coefficients $f_n$ and $g_n$. Then one can see using the transmission conditions ($u|_+=u|_-$ and $\p_\nu u|_+=k_j \p_\nu u |_-$) on the interfaces that
the solution to \eqnref{cond_eqn} is given by
\beq\label{sol_cond}
 \ds (u-h)(\xi,\Gt)=C+\begin{cases}
 \ds \sum_{n\neq 0} \left(a_n e^{|n| \xi} + b_n e^{-|n|\xi}\right)e^{in\theta} ,\quad &x\in\mathbb{R}^2 \setminus (D_1 \cup D_2),\\
 \ds \sum_{n\neq 0} \left(a_n e^{|n|\xi}+b_n e^{|n|(2\xi_0+\xi)}\right)e^{in\theta},& x\in D_1,\\
 \ds  \sum_{n\neq 0}\left(a_n e^{|n|(2\xi_0-\xi)}+b_n e^{-|n|\xi}\right)e^{in\theta},& x\in D_2,\\
 \end{cases}
\eeq
where
\beq\label{anbn}
a_n =\frac{g_n-\tau_1^{-1} e^{2|n|\xi_0}f_n }{\tau^{-1}e^{4|n|\xi_0}-1}, \quad
b_n= \frac{f_n-\tau_2^{-1} e^{2|n|\xi_0}g_n}{\tau^{-1} e^{4|n|\xi_0}-1},
\eeq
and
$$
C=-\sum_{n\ne 0} (a_n + b_n )
$$
Here, we set
\beq
\tau=\tau_1\tau_2,\quad \tau_j=\frac{k_j-1}{k_j+1},\quad j=1,2.
\eeq

One can see that the equation \eqnref{cond_eqn} is solvable if and only if the denominator in \eqnref{anbn} is non-zero, namely,
$$
\tau^{-1} e^{4n\xi_0}-1 \neq 0, \quad n=1,2, \ldots.
$$
This condition can be rewritten as
\beq\label{soltwo}
\Gl(k_1) \Gl(k_2) \neq \left(\pm \frac{1}{2} e^{-2n\xi_0} \right)^2, \quad n=1,2, \ldots.
\eeq
We emphasize that $\pm \frac{1}{2} e^{-2n\xi_0}$ are eigenvalues of the NP operator for two disks (see \cite{LiYu14}).

Note that
\beq
V_n:= \left\{ (k_1, k_2) \in \Cbb^2 ~|~ \Gl(k_1) \Gl(k_2) = \left(\frac{1}{2} e^{-2n\xi_0} \right)^2 \right\}, \quad n=1, 2, \ldots,
\eeq
is a two dimensional algebraic variety in $\Cbb^2$. So, the insolvability region in presence of two disks is
a countable family of two dimensional varieties in $\Cbb^2$.
It suggests that it is difficult to find a solvability condition for the free space problem or the boundary value problem when complex conductivities of inclusions are different. However, since $\Gl$ maps the right half plane onto the outside of the disk of radius $1/2$ (centered at $0$) and $(\frac{1}{2} e^{-2n\xi_0})^2 < 1/4$, we see that if $k_1' \ge 0$ and $k_2' \ge 0$ then $(k_1, k_2) \notin V_n$. This is the standard case of elliptic equations. On the other hand, $\Gl$ maps the real line onto the real line and $(\frac{1}{2} e^{-2n\xi_0})^2 >0$, and hence if $k_1''$ and $k_2''$ have the same signs, then $(k_1, k_2) \notin V_n$. In the next section we show the problem \eqnref{NBPmulti} in presence of multiple inclusions is well-posed if the imaginary parts of $k_j$ have the same signs.

\subsection{Uniqueness and existence}

The purpose of this section is to present a condition on $k_j$ under which the problem \eqnref{NBPmulti} is well-posed. The condition is that the imaginary parts of $k_j$'s are either all positive or all negative (they can be zero if the real part is non-negative). It is worth mentioning that if the real parts of $k_j$ are all non-negative, \eqnref{NBPmulti} is elliptic and well-posed.

Let
$$
S_L^+:= S_L \cap \{ k'+ik'' : k'\in \Rbb, \ k'' \ge 0\},
$$
where $S_L$ is the sector defined in subsection \ref{sec:uni}.
Let
\beq
\GT^+:=\bigcup_{L>0} \left(S_L^+\right)^M.
\eeq
It is worth mentioning that $\GT^+$ is the collection of all $\Bk =(k_1, \ldots, k_M) \in \Cbb^M$ with  $k_j\in \mathbb{C}\setminus (-\infty,0)$ and $k_j''\ge0$ (for all $j$). We define $\GT^-$ likewise, and let
\beq
\GT:= \GT^+ \cup \GT^-.
\eeq

We have the following well-posedness result for $\Bk \in \GT$.

\begin{theorem}\label{thm:invmulti}
If $\Bk=(k_1, \ldots, k_M) \in \GT$, then $\Dbb(\GL) - \Kbb^*$ is invertible on $\BH$, where $\GL=(\Gl_1, \ldots, \Gl_M)=(\Gl(k_1), \ldots, \Gl(k_M))$.
\end{theorem}
\pf
We realize $\Dbb(\GL) - \Kbb^*$ as a compact perturbation of invertible operator. In fact, we write
$$
\Dbb(\GL) - \Kbb^* = \mbox{diag} \left[ \Gl_1 I-\Kcal_{1}^*, \ldots, \Gl_M I- \Kcal^*_{M} \right]
- \begin{bmatrix} 0 & \p_{\nu_1} \Scal_{2} & \cdots & \p_{\nu_1} \Scal_{M}\\
\p_{\nu_2} \Scal_{1} & 0 & \cdots & \p_{\nu_2} \Scal_{M} \\
\vdots & \vdots & \ddots & \vdots\\
\p_{\nu_M} \Scal_{1} & \p_{\nu_M} \Scal_{2} & \cdots & 0 \end{bmatrix}.
$$
Since there is a positive distance between inclusions, the second operator on the righthand side above is a compact operator on $\BH$, while the first operator is invertible. So to show invertibility of $\Dbb(\GL) - \Kbb^*$ it suffices to prove its injectivity by the Fredholm alternative. Suppose that $\Phi=[\Gvf_1, \ldots, \Gvf_M]^T$ satisfies $(\Dbb(\GL) - \Kbb^*)[\Phi]=0$. Then we can see that the function $u$, defined by
$$
u(x) = \sum_{j=1}^M \Scal_{j} [\Gvf_j] (x), \quad x \in \Rbb^d,
$$
is a solution to
$$
\begin{cases}
\ds\nabla \cdot \Big( \sum_{j=1}^M k_j \chi(D_j) + \chi(\Rbb^d \setminus D_j) \Big) \nabla u =0 \quad\mbox{in } \Rbb^d, \\
\ds u(x) = O(|x|^{-d+1}) \quad\mbox{as } |x| \to \infty.
\end{cases}
$$
So, we have
\beq\label{total}
\int_{\Rbb^d \setminus \cup_{j=1}^M D_j} |\nabla u(x)|^2 dx + \sum_{k_j' \neq 0} k_j' \int_{D_j} |\nabla u(x)|^2 dx + i \sum_{j=1}^M k_j'' \int_{D_j} |\nabla u(x)|^2 dx=0.
\eeq
In particular, we have
\beq\label{imaginary}
\sum_{j=1}^M k_j'' \int_{D_j} |\nabla u(x)|^2 dx=0.
\eeq

Suppose that $\Bk \in \GT^+$. Since $k_j'' \ge 0$ for all $j$ and $k_j''>0$ if $k_j' < 0$, we infer from \eqnref{imaginary} that $u$ is constant in $D_j$ if $k_j' < 0$. Then \eqnref{total} becomes
\begin{align*}
\int_{\Rbb^d \setminus \cup_{j=1}^M D_j} |\nabla u(x)|^2 dx + \sum_{k_j' >0} k_j' \int_{D_j} |\nabla u(x)|^2 dx=0,
\end{align*}
from which we infer that $u$ is constant in $\Rbb^d \setminus \cup_{j=1}^M D_j$. Since $u(x) \to 0$ as $|x| \to \infty$, we have $u(x)=0$ in $\Rbb^d \setminus \cup_{j=1}^M D_j$. By the transmission condition (continuity of the potential) on $\p D_j$, we have $u=0$ in $D_j$, $j=1, \ldots, M$. Thus we have $\Gvf_j= \p_\nu u|_+ - \p_\nu u|_{-} = 0$ on $\p D_j$. This completes the proof.
\qed

The following is the main theorem of this subsection.

\begin{theorem}\label{thm:wellmulti}
For all $\Bk=(k_1, \ldots, k_M) \in \GT$ there is a unique solution $u_\Bk \in H^1(\GO)$ to \eqnref{NBPmulti}, and there is a unique $\psi^\Bk_j \in H^{-1/2}_0(\p D_j)$, $j=1, \ldots, M$, such that $u_\Bk$ is represented as
\beq\label{Urepmulti}
u_{\Bk}(x) = U(x) - \sum_{j=1}^M \Ncal_{j} [\psi^\Bk_j] (x), \quad x \in \GO,
\eeq
where $U$ is the solution to \eqnref{PDE}.
\end{theorem}
\pf
Uniqueness of the solution $u_\Bk$ can be proved similarly to Theorem \ref{thm:invmulti}. Existence and uniqueness of $\psi_j^\Bk$ can be proved in a similar way to the proof of Theorem \ref{thm:well} using Theorem \ref{thm:invmulti}.
\qed

\subsection{Uniformity of estimates and Lipschitz dependency}

We want to extend uniformity results in previous sections to the case of multiple inclusions. If $k_1, \ldots, k_M$ are identical, so are $\Gl_1, \ldots, \Gl_M$, say $\Gl$. Then, the equation \eqnref{solmatrix} takes the form
$$
\left(\Gl \BI - \Kbb^*\right)[\Phi_\Bk] = \p h_\Bk,
$$
where $\BI$ is the identity operator, and the spectral resolution can be applied as in section \ref{sec:uni}. However, if $\Gl_i \neq \Gl_j$ for some $i,j$, we do not have the spectral resolution and arguments of the previous sections cannot be applied. We overcome this difficulty by induction.

For $\Bk=(k_1, \ldots, k_M)$, let $\Gl(\Bk):= (\Gl(k_1), \ldots, \Gl(k_M))$ where $\Gl$ is the bilinear transform defined by \eqnref{bilinear}. For $L>0$ set $\GT_L:=(S_L^+)^M $ and
\beq
\GS_L:= \{ \, \GL= \Gl(\Bk) \, | \, \Bk \in \GT_L  \, \}.
\eeq

We first prove the following lemma.
\begin{lemma}\label{lemma3.1-1}
For any $R>0$ and $L>0$ there is a constant $C=C(L, R)$ such that
\beq\label{lemma3.1}
\left\| (\BI - \Dbb(\GL)^{-1} \Kbb^* )^{-1} \right\| \leq C,\quad \left\| (\BI - \Kbb^*\Dbb(\GL)^{-1}  )^{-1} \right\| \leq C
\eeq
for any $\GL=(\Gl_1, \ldots, \Gl_M) \in \GS_L$ with $|\GL| \leq R$. Here $\| \cdot \|$ is the operator norm on $\BH$.
\end{lemma}
\pf
For $\GL \in \GS_L$, $\Dbb(\GL) - \Kbb^*$ is invertible on $\BH$. So,
\beq
f(\GL): = \left\| (\BI - \Dbb(\GL)^{-1} \Kbb^* )^{-1} \right\|
\eeq
is well-defined. Since $\GS_L$ is a closed set,  it is sufficient to show that $f$ is continuous to obtain the first estimate in \eqnref{lemma3.1}.

Fix $\GL\in \GS_L$ and let $\{\GL_n\}$ be a sequence in $\GS_L$ converging to $\GL$.  Note that
\begin{align}
|f(\GL_n)-f(\GL)| &\le \left\| (\BI - \Dbb(\GL_n)^{-1} \Kbb^* )^{-1} - (\BI - \Dbb(\GL)^{-1} \Kbb^* )^{-1} \right\| \nonumber \\
&\leq  \left\| (\BI - \Dbb(\GL_n)^{-1} \Kbb^* )^{-1} \right\| \left\|(\Dbb(\GL)^{-1}-\Dbb(\GL_n)^{-1})\Kbb^* \right\| \left\| (\BI - \Dbb(\GL)^{-1} \Kbb^* )^{-1} \right\| \nonumber\\
&=f(\GL_n) \left\|(\Dbb(\GL)^{-1}-\Dbb(\GL_n)^{-1})\Kbb^* \right\| f(\GL). \label{diff_lambda}
\end{align}
Since $\|(\Dbb(\GL)^{-1}-\Dbb(\GL_n)^{-1})\Kbb^*\|\rightarrow 0$ as $n\rightarrow \infty$, we have
\begin{equation*}
\left\| (\Dbb(\GL)^{-1}-\Dbb(\GL_n)^{-1})\Kbb^* \right\| f(\GL) < \frac12
\end{equation*}
for all sufficiently large $n$, and hence we obtain from \eqref{diff_lambda} and the triangular inequality that
\begin{equation*}
f(\GL_n) < 2 f(\GL).
\end{equation*}
Inserting it into \eqref{diff_lambda}, it now follows
\begin{align*}
|f(\GL_n)-f(\GL)| \leq  2 \left\|(\Dbb(\GL)^{-1}-\Dbb(\GL_n)^{-1})\Kbb^* \right\| f(\GL)^2.
\end{align*}
Therefore $f$ is continuous on $\GS_L$. Similarly, we obtain the second estimate in \eqnref{lemma3.1}.
\qed

\begin{theorem}\label{thm_main}
For any $L>0$ there is a constant $C=C(L)$ such that
\begin{equation}\label{thm:bound2}
\|(\BI - \Dbb(\GL) ^{-1}\Kbb^* )^{-1}\| \leq C ,\quad \|(\BI - \Kbb^*\Dbb(\GL) ^{-1} )^{-1}\| \leq C
\end{equation}
for all $\GL \in \GS_L$.
\end{theorem}
\pf We only prove the first estimate in \eqnref{thm:bound2} because the second one can be proved similarly. We use  induction on $M$, the number of inclusions.

If $M=1$, \eqnref{specresol1} shows
\beq
(I - \Gl_1^{-1} \Kcal^*_1)^{-1} = \int_{\Gs(\Kcal_1^*)} \frac{1}{1- \Gl_1^{-1}t} d\Ecal_1(t)
\eeq
where $\Ecal_1(t)$ is the resolution of identity for $\Kcal_1^*$ on $\p D_1$. Since $|1- \Gl_1^{-1}t| \ge C$, we see that $\|( I - \Gl_1^{-1} \Kcal_{\p D_1}^*)^{-1}\| \le C$ uniformly in $\Gl_1 \in O_L$. Assume $M>1$. Without loss of  generality, we may assume that $|\lambda_1| \leq \cdots \leq |\lambda_M|$ by reordering $\lambda_j$'s. Let
$$
\GL_M:= (\Gl_1, \ldots, \Gl_M) \quad\mbox{and}\quad \GL_{M-1}:= (\Gl_1, \ldots, \Gl_{M-1}).
$$
Denote by $\Kbb_M^*$ the NP operator given in \eqnref{NPmulti} when there are $M$ inclusions. Then it takes the following form:
\begin{equation*}
\Kbb_M^* = \begin{bmatrix}   &  &  & \p_{\nu_1} \Scal_{M}\\
 & \Kbb_{M-1}^* &  & \p_{\nu_2} \Scal_{M} \\
 &  &  & \vdots\\
  \p_{\nu_M} \Scal_{1} & \p_{\nu_M} \Scal_{2} & \cdots & \Kcal^*_{M}\end{bmatrix} = :
  \begin{bmatrix}    \Kbb_{M-1}^* &  \Ebb  \\
    \Fbb &  \Kcal^*_{M}\end{bmatrix}.
\end{equation*}
We emphasize that $\Ebb$ is a bounded operator from $H^{-1/2}(\p D_M)$ into $H^{-1/2}(\p D_1) \times \ldots \times H^{-1/2}(\p D_{M-1})$ and $\Fbb$ is from $H^{-1/2}(\p D_1) \times \ldots \times H^{-1/2}(\p D_{M-1})$ into $H^{-1/2}(\p D_M)$.
Let
$$
\Phi_M=\begin{bmatrix} \Gvf_1 \\ \vdots \\ \Gvf_M \end{bmatrix}, \quad \Phi_{M-1}=\begin{bmatrix} \Gvf_1 \\ \vdots \\ \Gvf_{M-1} \end{bmatrix}, \quad \Psi_M=\begin{bmatrix} \psi_1 \\ \vdots \\ \psi_M \end{bmatrix}, \quad \Psi_{M-1}=\begin{bmatrix} \psi_1 \\ \vdots \\ \psi_{M-1} \end{bmatrix},
$$
and suppose that the following integral equation holds:
\begin{align}
\left( \BI - \Dbb(\GL_M)^{-1} \Kbb_M^* \right)[\Phi_M] =  \Psi_M.  \label{1}
\end{align}
Then we may rewrite this equation as
\begin{align}
\left( \BI-\Dbb(\GL_{M-1})^{-1} \Kbb_{M-1}^* \right)[\Phi_{M-1}] - \Dbb(\GL_{M-1})^{-1}\Ebb[\Gvf_M] & = \Psi_{M-1}, \label{2} \\
-\Gl_M^{-1}\Fbb[\Phi_{M-1}] + ( I  -\Gl_M^{-1} \Kcal^*_{M})[\Gvf_M] & = \psi_M. \nonumber
\end{align}
Solving the first equation for $\Phi_{M-1}$ and substituting it into the second equation yield
\begin{align*}
&(I - \Gl_M ^{-1} \Kcal^*_{M})[\Gvf_M] - \Gl_M ^{-1}\Fbb \, \left( \BI-\Dbb(\GL_{M-1})^{-1} \Kbb_{M-1}^* \right)^{-1}\Dbb(\GL_{M-1})^{-1}  \Ebb[\Gvf_M] \\
& =   \psi_M + \Gl_M ^{-1}\Fbb \, \left( \BI-\Dbb(\GL_{M-1})^{-1} \Kbb_{M-1}^* \right)^{-1} [\Psi_{M-1}].
\end{align*}
The induction hypothesis is
\beq\label{induction}
\left\| (I - \Gl_M ^{-1} \Kcal^*_{M})^{-1} \right\| \le C, \quad \left\| \left( \BI-\Dbb(\GL_{M-1})^{-1} \Kbb_{M-1}^* \right)^{-1} \right\| \le C.
\eeq
So, we infer that
\begin{align*}
\| \Gvf_M \| \le C_1 \Gl_M^{-1} \| \Gvf_M\| + C_2 \|\psi_M\| + C_3 \lambda_M^{-1} \| \Psi_{M-1} \|
\end{align*}
for some constants $C_1$, $C_2$, and $C_3$ independent of $\GL$. Here the norms $\| \cdot \|$ are $H^{-1/2}$-norms on appropriate boundaries $\p D_j$. For example, $\| \Gvf_M \|$ is on $\p D_M$ and $\| \Psi_{M-1} \|$ is on $\p D_1 \times \cdots \times \p D_{M-1}$. So, if $|\Gl_M| \ge R$ for some large $R$, then we have
$$
\| \Gvf_M \| \le C_4 \| \Psi_M \|.
$$
We then see from \eqnref{2} and \eqnref{induction} that
$$
\| \Phi_M \| \le C \| \Psi_M \|
$$
for some constant $C$. So we obtain
$$
\|(\BI - \Dbb(\lambda) ^{-1}\Kbb^* )^{-1}\| \leq C
$$
when $|\Gl_M| \ge R$. Combined with Lemma \ref{lemma3.1-1}, the proof is completed.  \qed

As a consequence we have the following corollary.
\begin{cor}\label{thm3.3}
For each $L>0$ there is a constant $C=C(L)$ such that
\begin{equation}\label{thm1}
\|(\Dbb(\Gl(\Bk)) - \Kbb^* )^{-1} - (\Dbb(\Gl(\Bs))  - \Kbb^* )^{-1}\| \leq C\sum_{j=1}^M \frac{|k_j - s_j|}{(1+ |k_j|)(1+ |s_j|)}
\end{equation}
for all $\Bk, \Bs \in \GT_L$.
\end{cor}
\pf Using \eqnref{thm:bound2} and the identity
\begin{align*}
&(\Dbb(\Gl(\Bk))  - \Kbb^* )^{-1} - (\Dbb(\Gl(\Bs))  - \Kbb^* )^{-1}\\
&= (\BI - \Dbb(\Gl(\Bk)) ^{-1}\Kbb^* )^{-1} (\Dbb(\Gl(\Bk)) ^{-1}-\Dbb(\Gl(\Bs)) ^{-1})(\BI - \Kbb^*\Dbb(\Gl(\Bs)) ^{-1} )^{-1},
\end{align*}
we have
\begin{align*}
\|(\Dbb(\Gl(\Bk))  - \Kbb^* )^{-1} - (\Dbb(\Gl(\Bs))  - \Kbb^* )^{-1}\|
 &\leq C \|\Dbb(\Gl(\Bk)) ^{-1}-\Dbb(\Gl(\Bs)) ^{-1}\| \\
&\leq C\sum_{j=1}^M \frac{|k_j - s_j|}{(1+ |k_j|)(1+ |s_j|)}\end{align*}
for all $\Bk, \Bs \in \GT_L$.
\qed

Now we obtain the uniformity of the regularity estimate and the Lipschitz continuity estimate for the solution $u_{\Bk}$ to \eqnref{NBP} when there are multiple inclusions. The following theorem can be proved in a similar way to the proof of Theorem \ref{thm:uni} and Theorem \ref{conti2} using \eqnref{thm:bound2} and \eqnref{thm1}.

\begin{theorem}\label{thm:uni3}
For any $L>0$ there is a constant $C=C(L)$ such that
\beq\label{Uni}
\|u_\Bk\|_{H^1(\GO)} \leq C \|g\|_{-1/2(\p \GO)}
\eeq
and
\begin{equation}\label{Lip}
\|u_\Bk- u_\Bs\| _{H^1(\GO)} \leq C\sum_{j=1}^M \frac{|k_j - s_j|}{(1+ |k_j|)(1+ |s_j|)} \|g\|_{-1/2(\p \GO)}
\end{equation}
for all $\Bk, \Bs \in \GT_L$ and $g \in H^{-1/2}_0(\p \GO)$.
\end{theorem}

The results of this section can also be applied to extend Theorem \ref{uniitrans} and the uniformity of the boundary perturbation formula (Theorem \ref{uniform}) to the multiple inclusion case. But we do not state the extended results. We only mention that there are two different boundary perturbation formulas: one for the case of well-separated diametrically small inclusions and another for the case of a cluster of closely located small inclusions (see \cite{AKKL05}). The uniform validity of the boundary perturbation formula for the first case can be proved using arguments in section \ref{sec:asymp}, and that for the second case can be proved using results of this section.


\end{document}